\documentclass[11pt]{amsart}
\usepackage{amssymb} 
\usepackage{amsmath} 
\usepackage{graphicx} 
\usepackage[english]{babel} 
\usepackage{epsfig}
\usepackage{enumerate}
\usepackage{color}
\usepackage{mathtools}
\usepackage[all]{xy}
\usepackage{mathrsfs}
\usepackage{graphics}
\usepackage{amsthm}
\numberwithin{equation}{subsection}

\usepackage[usenames,dvipsnames]{pstricks}
 \usepackage{pst-grad} 
\usepackage{pst-plot} 

\swapnumbers

\newenvironment{proo}{\begin{trivlist} \item{\emph{Proof.}}}
 {\hfill $\square$ \end{trivlist}}
\swapnumbers

\theoremstyle{definition}
\newtheorem{definition}[subsubsection]{Definition}
\newtheorem{example}[subsubsection]{Example}
\newtheorem{notation}[subsubsection]{Notation}
\newtheorem{remark}[subsubsection]{Remark}
\theoremstyle{plain}
\newtheorem{theorem}[subsubsection]{Theorem}
\newtheorem{proposition}[subsubsection]{Proposition}
\newtheorem{lemma}[subsubsection]{Lemma}
\newtheorem{corollary}[subsubsection]{Corollary}

\def\Ordm{{\mbox {Ord}^m}}

\def\ZZ{{\mathbb{Z}}}
\def\K{{\mathbb K}}
\def\ot{\otimes}

\def\noi{\noindent}
\def\Dy{\mbox {\it Dyck}}
\def\Dyc{\mbox {\it Dy}}
\def\D{\mathcal{D}}

\def\ov{\overline}

\def\up{\mathcal{UP}}
\def\dw{\mathcal{DW}}
\def\t{\times}
\def\ot{\otimes}
\def\D{\mathcal D}
\def\Dot2{\ov{\D_m^+\ot \D_m^+}}

\def\lam{{\underline {\lambda}}}
\def\Dot2{\ov{\D_m^+\ot \D_m^+}}

\begin{document}

\author[D. L\'opez N., L.-F. Pr\'eville-Ratelle, M. Ronco]{Daniel L\'opez N., Louis-Fran\c cois Pr\'eville-Ratelle, Mar\'\i a Ronco}
\address{DL: Institut de Math\'ematiques de Jussieu-Paris Rive Gauche\\
B\^atiment Sophie Germain, 8 Place Aur\'elie Nemours, 75205 PARIS Cedex 13\\France}
\email{daniel.lopez@imj-prg.fr}
\address{LFPR: Instituto de Matem\'atica y F\'\i sica\\
Campus Norte, Camino Lircay s/n
\\Talca, Chile}
\email{lfprevilleratelle@inst-mat.utalca.cl}
\address{MOR: IMAFI, Universidad de Talca\\ Campus Norte, Avda. Lircay s/n\\ Talca, Chile}
\email{maria.ronco@inst-mat.utalca.cl}

\title{A simplicial complex spliting associativity}


{\thanks{Our joint work was partially supported by the projects Fondecyt Postdoctorado 3140298, Fondecyt Regular 1171209 and MathAmSud 17-Math-05 LIETS.
}}


\begin{abstract} We introduce a simplicial object $(\{ \Dy^m\}_{m\geq 0}, {\mathbb F}_i, {\mathbb S}_j)$ in the category of non-symmetric algebraic operads, satisfying that $\Dy^0$ is the operad of associative algebras and $\Dy^1$ is J.-L. Loday\rq s 
operad of dendriform algebras. The dimensions of the operad $\Dy^m$ are given by the Fuss-Catalan numbers.

Given a family of partially ordered sets ${\bold P}=\{P_n\}_{n\geq 1}$ we show that, under certain conditions, the vector space spanned by the set of $m$-simpleces of ${\bold P}$ is a $\Dy^m$ algebra. This construction, applied to certain combinatorial Hopf algebras, whose associative product  comes from a dendriform structure,  provides examples of $\Dy^m$ algebras. \end{abstract}

\maketitle

\section*{Introduction} \label{section:introduction} 
\medskip

In \cite{Lod}, J.-L. Loday introduced the notion of dendriform algebra as a type of associative algebra, whose product splits as the sum of two binary operations. Many associative algebras, as the algebras defined by shuffles (see \cite{EilMac}, \cite{LodRon}) and the Rota-Baxter algebras (see \cite{Agu}), are examples of dendriform algebras. 

Recent publications (see for instance \cite{GubKol}, \cite{Bre}, \cite{BaiBelGuoNi}, \cite{BaiGuoPei} and \cite{Gir}) deal with different ways of spliting associativity, involving the dendriform operad. 
We introduce a family of non-symmetric operads $\{\Dy^m\}_{m\geq 0}$, equipped with operad morphisms ${\mathbb F}_i: \Dy^m\longrightarrow \Dy^{m-1}$ and ${\mathbb S}_i:\Dy^m\longrightarrow \Dy^{m+1}$, for $0\leq i\leq m$, satisfying that: \begin{enumerate}
\item $\Dy^0$ is the operad ${\mbox{As}}$ of associative algebras and $\Dy^1$ is the operad of dendriform algebras,
\item  $(\{ \Dy^m\}_{m\geq 0}, {\mathbb F}_i, {\mathbb S}_j)$ is a simplicial object in the category of non-symmetric operads,
\item the dimension of the subspace of homogeneous elements of degree $n$ of the operad $\Dy ^m$ is the Fuss-Catalan number $d_{m,n}$.\end{enumerate}

As $\Dy^m$ is a non-symmetric operad, it is completely described by its free object over one element (see \cite{MaScSt} and \cite{LodVal}).  

For $m=1$, the Fuss-Catalan number $d_{1,n}$ coincides with the Catalan\rq s number $c_n= \frac 1{n+1}\binom{2n}n$, which is the cardinal of the set ${\mathcal Y}_n$ of planar binary rooted trees with $n{\mbox{+}1}$ leaves. 
In \cite{LodRon}, the free dendriform algebra over one element was described on the vector space $\K[{\mathcal Y}]:=\bigoplus_{n\geq 1}\K[{\mathcal Y}_n]$, spanned by the set of planar binary rooted trees. In \cite{LodRon1}, J.-L. Loday and the third author proved that the dendriform structure of $\K[{\mathcal Y}]$, is completely determined by the Tamari order $\leq_{Ta}$ (see \cite{Tam}), and two morphisms of partially ordered sets $/, \backslash: {\mathcal Y}_n\times {\mathcal Y}_r\longrightarrow {\mathcal Y}_{n+r}$, for $n,r\geq 1$.  

The last result also holds for other dendriform algebras, whose underlying vector spaces admit a graded basis $\bigcup_{n\geq 1}P_n$, where $P_n$ is a partially ordered set, for $n\geq 1$. Let us mention the dendriform algebras spanned by\begin{enumerate}
\item the sets $\Sigma _n$ of permutations of $n$ elements equipped with the weak Bruhat order (see \cite{LodRon1}), 
\item the sets of surjective maps ${\mbox{Surj}_n}$, equipped with the facial order introduced in \cite{KrLaNoFhSC}. Dendriform structures on the vector space spanned by the set $\bigcup _{n\geq 1}{\mbox{Surj}_n}$, were defined in \cite{Cha0} and in \cite{LodRon2}. 
\item the sets of planar rooted trees ${\mathcal T}_n$ , equipped with the partial order defined in \cite{PalRon}, which extends the Tamari order.\end{enumerate}

The functor ${\mbox {Simp}(P)}$ associates to any partially ordered set $(P,\leq)$ a simplicial set ${\mbox {Simp}(P)}$, whose $m$-simpleces are the $m$-tuples $(p_1\leq \dots \leq p_m)$ in $P^{m}$.  Given a family ${\bold P}=\{ P_n\}_{n\geq 1}$ of partially ordered sets and $m\geq 1$, we denote by ${\mbox {Simp}({\bold P})^m}$ the graded set of $m$-simpleces, whose elements of degree $n$ are the $m$-simpleces of $P_n$, for $n\geq 1$. In Section 3, we introduce the notion of dendriform poset. The vector space spanned by a dendriform poset is a dendriform algebra, while the vector space spanned by the set ${\mbox {Simp}({\bold P})^m}$, of $m$-simpleces of $\{\bold P\}$, has a natural structure of  $\Dy^m$ algebra. This result provides examples of $\Dy^m$ algebras. For instance, we prove that the families of partially ordered sets $\{{\mathcal Y}_n\}, \{ \Sigma _n\}, \{{\mbox{Surj}_n}\}$ and $\{{\mathcal T}_n\}$ described in the paragraph above, are dendriform posets. 

On the other hand, motivated by the combinatorics of the Garsia-Haiman spaces, F. Bergeron introduced in \cite{Ber} a generalization of the Tamari lattice, called the $m$-Tamari lattice, defined on the sets of $m$-Dyck paths. In the last section, we describe a $\Dy^m$ algebra structure on the vector space spanned by the set $\bigcup_{n\geq 1}\Dyc_n^m$, of all $m$-Dyck paths. We prove that the products $*_0,\dots ,*_m$, which define the $\Dy^m$ algebra, are described by intervals of the $m$-Tamari order.

In \cite{NovThi}, J.-C. Novelli and J.-Y. Thibon introduced the combinatorial Hopf algebra $^m{\mbox{\bf{FQSym}}}$ of $m$-permutations. Their construction provide $m$-analogues of a large family of combinatorial Hopf algebras, in particular they defined a combinatorial Hopf algebra $^m{\mbox{\bf{PBT}}}$, whose underlying vector space is spanned by the set of planar rooted $m$-ary trees, and proved that its associative product is described by the $m$-Tamari lattice. J.-C. Novelli introduced in \cite{Nov} a family of operads, called $m$-dendriform operads, whose dimensions are the Fuss-Catalan numbers. Though the operad of $m$-dendriform algebras and $\Dy^m$ have the same dimensions, both operads are not isomorphic for $m\geq 2$. Clearly, using any bijective map between the set of $m$-Dyck paths and the set of $m$-ary planar rooted trees, we may define a $\Dy^m$ algebra structure on the vector space spanned by the second set. But the description of this $\Dy^m$ algebra is much more complicated than the one defined on $m$-Dyck paths.

Our motivation to work on this type of algebraic operads is two fold. \begin{enumerate}
\item Let ${\mathcal P}$ be an operad such that there exists, for $n\geq 1$, a partial order defined on a basis $P_n$ of the space ${\mathcal P}_n$, satisfying that these orders are compatible with the operad structure. 
An interesting problem is to study the existence of an operad $\Ordm({\mathcal P})$, spanned by the operations ${\mbox {Simp}({\bold P})^m}$, for $m\geq 1$, which must be a quotient of the Hadamard product ${\mathcal P}^{\otimes_Hm}$. 
The dendriform operad ${\mbox {Dend}}$ is an example of this type of operad, the set ${\mathcal Y}_n$ of planar rooted binary trees is a basis of ${\mbox {Dend}}_n$, and the Tamari order is compatible with the operad structure. Nevertheless, 
${\mbox{dim}_{\K}(\Dy_3^2)}=12$ and ${\mbox{dim}_{\K}({\mbox {Ord}^2}({\mbox{Dend}})_3)}$ is $13$, so even if the operad $\Ordm({\mbox{Dend}})$ exists, it is not $\Dy^m$. In recent publications, as \cite{FoiMal}, \cite{ChaPilPon}, \cite{Pil} and \cite{PilPon},
algebraic structures related to families of partially ordered sets were described, some of them give new examples of dendriform algebras.
\item in \cite{Cha}, F. Chapoton introduced a differential non-symmetric operad ${\mathcal K}$, whose algebras are dendriform algebras equipped with an additional associative product satisfying certain relations. 
The operad ${\mathcal K}$ is Hopf, which implies that any free ${\mathcal K}$ algebra has a natural structure of conilpotent bialgebra. As proved in \cite{BurRon}, the subspace of primitive elements of any ${\mathcal K}$ bialgebra is a ${\mathcal S}_2$ algebra, where ${\mathcal S}_2$  is the second filtration stage of the surjection operad ${\mathcal S}$ (see \cite{BerFre}), so the operad ${\mathcal S}_2$ is an $E_2$ operad. Moreover, there exists an equivalence between the category of conilpotent ${\mathcal K}$ bialgebras and the category of ${\mathcal S}_2$-algebras, which implies that the operad ${\mathcal S}_2$ is completely described in terms of the operad ${\mathcal K}$ (see \cite{BurRon}). In \cite{BatBer}, M. Batanin and C. Berger introduced 
the filtered colored operad ${\mathcal L}$, called the lattice path operad, whose $n^{th}$ filtration stage ${\mathcal L}_n$ is an $E_n$ operad, for $n\geq 1$. Our aim is to study non-symmetric Hopf operads ${\mathcal P}$ defined on the space spanned 
by all lattice paths, applying the ideas developed by X. Viennot and the second author in \cite{PRvie}. As the operad ${\mathcal L}$ is symmetric, we are interested in the properties of the operads defined on the subspaces of primitive elements of 
${\mathcal P}$ bialgebras and in their relationship with the operad ${\mathcal L}$.\end{enumerate}

\vfill\eject

{\it Contents} 
\medskip

The first section contains the definition of $\Dy^m$ algebras and the description of the free $\Dy^m$ algebra over one element. We prove that the dimensions of $\Dy^m$ are given by the Fuss-Catalan numbers, and compare the operad 
$\Dy^2$ to other operads having the same dimensions, in particular with the operad of $m$-dendriform algebras introduced in \cite{Nov}.

The simplicial object $(\{ \Dy^m\}_{m\geq 0}, {\mathbb F}_i, {\mathbb S}_j\})$ is described in Section 2. We show that the degeneracy functors ${\mathbb S}_j:{\mbox{}\Dy ^{m+1}{\mbox{-alg}}}\longrightarrow  {\mbox{}\Dy^m{\mbox{-alg}}}$ preserve free objects. 

In Section 3 we introduce the definition of dendriform poset, in such a way that the vector space spanned by any dendriform poset is a dendriform algebra. We prove that, for any dendriform poset ${\bold P}$, the vector space $\Ordm({\mbox{\bf P}})$ spanned by the $m$-simpleces of the partially ordered sets $P_n$, is a $\Dy^m$ algebra. Using previous results, we describe examples of $\Dy^m$ algebras coming from dendriform posets. 

We recall basic definitions and constructions of $m$-Dyck paths, and define binary products $*_0, \dots ,*_m$ on the space $\K[\Dyc^m]$, spanned by the set of $m$-Dyck paths, in Section 4. We prove that the operad $\Dy^m$ 
is entirely described by the data $(\K[\Dyc^m], *_0,\dots ,*_m)$, so the combinatorial properties of $m$-Dyck paths define completely the operad. Finally, we show the $m+1$ binary operations $*_0,\dots , *_m$
are given by intervals of F. Bergeron\rq s $m$-Tamari lattice.  
\medskip

\subsection*{Acknowledgements}  D. L\'opez N. and M. Ronco want to thank specially Prof. Antonio Laface for his interest and support. The second author would like to thank Luc Lapointe for many fruitful discussions.
\section*{Preliminaries}

All the vector spaces considered in the present work are over $\K$, where $\K$ is a field. For any set $X$, we denote by $\K[X]$ the vector space spanned by $X$. For any positive integer $n\geq 1$, we denote the set $\{1, \dots ,n\}$ by $[n]$.

\section{$\Dy^m$ algebras}
\medskip 

\subsection{Definition and basic properties}
\medskip

Dendriform algebras were introduced by J.-L. Loday in \cite{Lod}. However, the first example of this type of algebra appeared many years before in \cite{EilMac}, S. Eilenberg and S. MacLane introduced the shuffle product on simplicial complexes, using the Alexander-Whitney map, and defined a {\it half product}, denoted $\uparrow$, which splits it. In \cite{Agu}, M. Aguiar proved that Rota-Baxter algebras have a natural structure of dendriform algebras.

\begin{definition} \label{defdend} (see \cite{Lod}) A {\it dendriform algebra} over $\K$ is a vector space $A$ equipped with binary operations $\succ $ and $\prec $ satisfying the following conditions\begin{enumerate}
\item $x\succ (y\succ z) = (x\succ y + x\prec y) \succ z$,
\item $x\succ (y\prec z) = (x\succ y)\prec z $, 
\item $x\prec (y\succ z + y \prec z)=(x\prec y)\prec z$,\end{enumerate}
for $x,y,z \in A$.\end{definition}
\medskip

\begin{remark} \label{splitassoci} Suppose that we split $\succ$ as the sum of two binary operations $*_0$ and $*_1$. In order to describe new relations, in such a way that 
conditions $(1)$ and $(2)$ of Definiton \ref{defdend} are satisfied, we set\begin{enumerate}
\item $x*_0(y*_0 z) = (x*y)*_0z$,
\item $x*_0(y*_1 z)= (x*_0y)*_1z$,
\item $x*_1(y\succ z)=(x*_1y + x\prec z)*_1 y$,
\item $x*_i(y\prec z)=(x*_iy)\prec z$, for $i= 0,1$. \end{enumerate}

Note that the first three relations split relation $(1)$ of Definition \ref{defdend}, while the last relation splits the second one.

It is not difficult to verify that, if we perform an analogous procedure with $\prec$, that is we write $\prec=\circ_1+\circ_2$ and splits the relations of dendriform, we get the same type of algebra, with $3$ binary products.\end{remark}

Remark \ref{splitassoci} motivates the following definition. 

\begin{definition} \label{defDyckmalgebra} For $m\geq 1$, a {\it $\Dy ^m$ algebra} over $\K$ is a vector space $D$ equipped with $m+1$ binary operations 
$*_i: D\ot D\longrightarrow D$, for $0\leq i\leq m$, satisfying the following relations
\begin{equation}\label{eq1} x*_i(y*_j z) = (x*_i y)*_j z,\ {\rm for}\ 0\leq i < j\leq m,\end{equation}
\begin{equation}\label{eq2}\sum _{j=0}^i x*_i (y*_j z) = \sum _{k=i}^m (x*_k y)*_i z,\ {\rm for}\ 0\leq i\leq m,\end{equation}
for any elements $x, y $ and $z$ in $D$.\end{definition}
\medskip

Let ${\mbox {As}}$ denotes the operad of associative algebras. It is immediate to see that $\Dy^0={\mbox {As}}$. A $\Dy^1$ algebra is a dendriform algebra, for $\succ:= *_0$ and $\prec:=*_1$. 

Any $\Dy^m$ algebra $(D, *_0,\dots ,*_m)$ is an associative algebra,  with the product $*:=*_0+\dots +*_m$. 
It is immediate to verify that the vector space $D$, with the products $\succ ^k := *_0+\dots +*_k$ and 
$\prec ^k :=  *_{k+1}+\dots +*_m$ is a dendriform algebra, for any $-1\leq k\leq m$.
\medskip

Dendriform algebras cannot have units, see for instance \cite{Ron}. A straightforward calculation proves that this result is also true for $\Dy^m$ algebras in general.

\begin{remark} \label{other basis} Let $(D, *_0,\dots ,*_m)$ be a $\Dy^m$ algebra. Define operations $\circ_i:=*_0+\dots +*_i$, for $0\leq i\leq m$.

It is easy to verify that conditions \ref{eq1} and \ref{eq2} are equivalent to \begin{enumerate}
\item for $0\leq i<j\leq m$, 
\begin{equation*} x\circ _i(y\circ _j z) - (x\circ _i y)\circ _jz = x\circ_i(y\circ_{j-1}z) - (x\circ_i y)\circ_{j-1}z,\end{equation*}
\item $x\circ_0(y\circ_0 z)= (x\circ_m y)\circ_0 z$,
\item for $1\leq i\leq m$,
\begin{equation*} x\circ_i(y\circ _i z) = (x\circ_m y)\circ_i z - (x\circ_m y)\circ _{i-1} z +x\circ_{i-1}(y\circ_{i-1}z),\end{equation*}
\end{enumerate}
for $x, y$ and $z$ in $D$.\end{remark}

\subsection{The free $\Dy ^m$ algebra} 
\medskip

\begin{remark} \label{freeoverone} Let $\Dy^m=\bigoplus_{n\geq 1}{\Dy_n^m}$ denote the free $\Dy ^m$ algebra over one generator. As the operad $\Dy ^m$ is non-symmetric, the free $\Dy^m$ algebra over a vector space $V$ (see \cite{MaScSt} and \cite{LodVal}) is the vector space 
\begin{equation}\label{eq3} \Dy^m(V):=\bigoplus_{n\geq 1} {\Dy_n^m}\otimes V^{\otimes n},\end{equation} 
with the products 
\begin{equation*} t\otimes (v_1\otimes \dots \otimes v_n)\ *_i\ w\otimes (u_1\otimes \dots \otimes u_p)\ :=\ (t*_iw)\otimes (v_1\otimes \dots \otimes v_n\otimes u_1\otimes \dots \otimes u_m),\end{equation*}
for $t\in \Dy_n^m$, $w\in \Dy_p^m$, $v_1\otimes \dots \otimes v_n\in V^{\otimes n}$ and $u_1\otimes \dots \otimes u_p\in V^{\otimes p}$.\end{remark}

\begin{remark}\label{rem:freeDyckwithbinarytrees}  In Definition \ref{defDyckmalgebra} the equation \ref{eq2} may be written as
\begin{equation*} (x*_i y)*_iz = \sum _{j = 0}^i x*_i (y*_j z)-\sum _{j = i+1}^m (x*_j y)*_i z,\ {\rm for}\ 0\leq i\leq m,\end{equation*}
for $x, y$ and $z$ in $D$.\end{remark}

\begin{notation} \label{notbasictrees} For $n\geq 1$, let ${\mathcal Y}_{n-1}^m$ be the set of all planar binary rooted trees with $n$ leaves and the internal vertices colored by the 
elements of $\{ 0,\dots ,m\}$. The unique element of ${\mathcal Y}_{0}^m$ is the tree $\vert$, with one leaf and no internal vertex. 

Given two colored trees, $t$ and $w$, and an integer $0\leq i\leq m$, we denote by $t\vee _{i}w$ the colored tree obtained by connecting the roots of $t$ and $w$ to a new root colored by $i$, where $t$ is on the left side 
and $w$ on the right side of $t\vee_i w$.

For any internal vertex $v$ of a colored planar binary rooted tree $t\in {\mathcal Y}_{n-1}^m$, we denote by $t_v$ the colored subtree of $t$ whose root is $v$. \end{notation}

Note that any tree $t\in {\mathcal Y}_{n-1}^m$ may be written as $t=t^l\vee _it^r$, for a unique integer $0\leq i\leq m$ and unique trees $t^l\in {\mathcal Y}_{n_1}^m$, $t^r\in {\mathcal Y}_{n_2}^m$ such that $n_1+n_2=n{\mbox{-} 2}$.

\begin{definition} \label{def:basisBm} For $n\geq 2$, let ${\mathcal B}_n^{m}$ be the subset of all the elements $t$ in ${\mathcal Y}_{n-1}^m$ such that for any any subtree $t_v=t_v^l \vee _{i} t_v^r$ of $t$, the color of the root of $t_v^l$ is $j$, for some $ j > i$.\end{definition}
\medskip

For $n=1$, ${\mathcal B}_1^m= {\mathcal Y}_{0}^m$ is the set whose unique element is the tree $\vert$. 

Let ${\mathcal B}^m$ be the graded set $\bigcup _{n\geq 1}{\mathcal B}_n^m$. For any $t = t^l\vee _{i}t^r\in {\mathcal B}^m$, the trees $t^l$ and $t^r$ belong to ${\mathcal B}^m$.
\medskip

\begin{definition} \label{prodsontrees} Let ${\frak D}^m$ be the graded vector space whose basis is the set $\bigcup _{n\geq 1}{\mathcal B}_n^m$. 

For any pair of trees $t\in {\mathcal B}_n^m$ and $w\in {\mathcal B}_r^m$, and any integer $0\leq i\leq m$, the product $t*_i w$ is defined recursively as follows,\begin{enumerate}
\item for ${\mbox{}n = r = 1}$, we set $\vert *_i\vert :=\vert \vee_i\vert$.
\item for $t = t^l\vee_{j} t^r$, or $n = 1$, with $i < j\leq m$, define $t*_i w := t\vee _{i}w \in {\mathcal B}_{n+r}^m$,
\item suppose that the products $t\rq *_iw\rq$ are defined for any pair of trees $t\rq \in {\mathcal B}_{n_1}^m$ and $w\rq \in {\mathcal B}_{r_1}^m$, for $n_1+r_1<n+r$. Moreover, we may also assume that, if 
$t\rq *_iw\rq= \sum _u \alpha_{t\rq w\rq}^u u$, then, for any $u\in {\mathcal B}_{n_1+r_1}^m$ such that $\alpha_{t\rq w\rq}^u\neq 0$, the color of the root of $u$ is largest or equal than the minimal element of the set $\{ i, s\}$, where $s$ is the color of the root of $t\rq$.

Let $t = t^l\vee_{j} t^r$, with $0\leq j\leq i$.\begin{enumerate}
\item for $j < i$, define $t *_i w := t^l\vee_{j}(t^r*_{i} w)$.
\item for $j = i$, applying Remark \ref{rem:freeDyckwithbinarytrees}, we get 
\begin{equation*} t*_iw = \sum _{k = 0}^i t^l *_i (t^r *_k w)- \sum _{k=i+1}^m ( t^l *_k t^r)*_i w.\end{equation*} 
Moreover, 

\noindent $(i)$ as $t\in {\mathcal B}_n^m$, the color $s$ of the root of $t^l$ satisfies that $i< s$. So, the element $\sum _{k = 0}^i t^l \vee_i (t^r *_k w)$ belongs to ${\mathfrak D}^m$,

\noindent $(ii)$ the product $t^l*_kt^r$ is a sum of trees, satisfying that the colors of their roots are largest or equal to ${\mbox {min}\{ k, s\}}$, where $s$ is the color of the root of $t^l$. As both integers are largest than $i$, the roots of all trees appearing 
in $t^l*_kt^r$ are largest than $i$.

So, the element $t*_iw$ is defined as
\begin{equation*}\qquad t*_iw := \sum _{k = 0}^i t^l \vee_i (t^r *_k w)- \sum _{k=i+1}^m ( t^l *_k t^r)\vee_i w.\end{equation*}
 \end{enumerate}\end{enumerate}
 It is easily seen that the products $*_i$ are defined on ${\mathfrak D}^m$ in such a way that $({\mathfrak D}^m, *_0,\dots ,*_m)$ is a $\Dy^m$ algebra.\end{definition}

\begin{example} \label{prodintrees} Let 

\begin{figure}[h]
\includegraphics[scale=0.6]{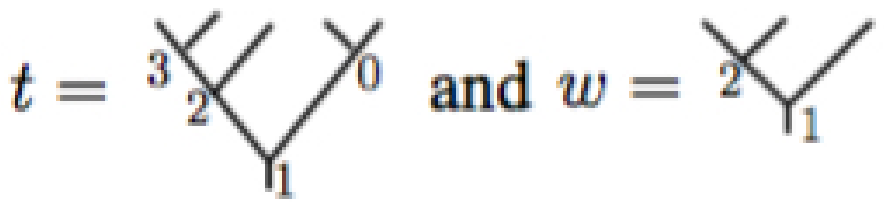}
\end{figure}

we get that 

\begin{figure}[h]
\includegraphics[scale=0.6]{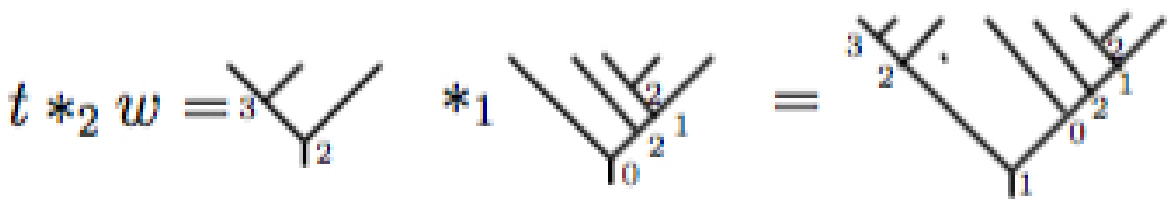}
\end{figure}

\end{example}
\medskip

The following result is a consequence of Remark \ref{rem:freeDyckwithbinarytrees} and Definition \ref{prodsontrees}.

\begin{proposition}\label{prop:freeDyckmontrees} The graded vector space ${\mathfrak D}^m$, equipped with the products $*_0,\dots ,*_m$, is the free $\Dy ^m$ algebra over one element $\#$.\end{proposition}
\medskip

\begin{proo} The linear homomorphism $\K\# \longrightarrow {\mathfrak D}^m$, maps $\#\mapsto \vert$.

Given a $\Dy^m$ algebra $(D, *_0^D,\dots ,*_m^D)$, a linear map $f$ from $\K\#$ to $D$ is completely determined by the element $f(\#)\in D$. There exists a unique linear homomorphism ${\tilde f}: {\mathfrak D}^m\longrightarrow D$ satisfying that\begin{enumerate}[(a)]
\item ${\tilde f}(\vert ):=f(\#)$,
\item for $t=t^l\vee _it^r\in {\mathcal B}_n^m$, 
\begin{equation*} {\tilde f}(t) = {\tilde f}(t^l)*_i^D{\tilde f}(t^r).\end{equation*}\end{enumerate}
As the elements of the set $\bigcup _{n\geq 1}{\mathcal B}_n^m$ do not contain any subtree of the form $(t\vee _i w)\vee_j u$, for $0\leq i\leq j \leq m$, we have that ${\tilde f}$ is well-defined. On the other hand, as $D$ is a $\Dy^m$ algebra,
the way products $*_i$ are defined, implies that ${\tilde f}$ is a $\Dy^m$ algebra homomorphism.
\end{proo}
\medskip

For any vector space $V$, equation \ref{eq3} shows that a basis of the free algebra ${\frak D}^m(V)$ over $V$ is given by
\begin{equation*} {\mathcal B}^m(V):= \bigcup_{n\geq 1} {\mathcal B}_n^m\times X^n,\end{equation*}
where $X$ is a basis of $V$.

\subsection{The generating series associated to $\Dy^m$} 
\medskip

Let 
\begin{equation}\label{eq4} d_{m,n}:= \frac 1{mn+1}\binom {(m+1)n}n,\end{equation} be the Fuss-Catalan number, for $m,n\geq 1$. 

The integer $d_{m,n}$ is the number of planar rooted $m$-ary trees, or equivalently the number of $m$-Dyck paths of length $n$, see for instance \cite{Ber}, \cite{BouFusPre}, \cite{NovThi} or \cite{Gir}. The generating series of $d_{m,n}$ is
\begin{equation*} x\cdot d_m(x)^{m+1} = d_m(x)-1.\end{equation*}
We want to see that the number of elements of ${\mathcal B}_n^{m}$ is $d_{m,n}$, for $m,n\geq 1$. 
\medskip

\begin{notation} \label{not:combs} Given a family of colored trees $t_1,\dots ,t_p$ and a family of integers $0\leq i_1,\dots , i_p\leq m$, we denote by \begin{enumerate}
\item $\Omega_{i_1,\dots ,i_p}^L(t_1,\dots ,t_p)$ the colored tree 
\begin{equation*}\Omega_{i_1,\dots ,i_p}^L(t_1,\dots ,t_p) := (((\vert \vee_{i_{p}} t_p)\vee_{i_{p-1}} t_{p-1})\dots )\vee_{i_1} t_1,\end{equation*}
\item $\Omega_{i_1,\dots ,i_p}^R(t_1,\dots ,t_p)$ the colored tree 
\begin{equation*}\Omega_{i_1,\dots ,i_p}^R(t_1,\dots ,t_p) := t_1\vee _{i_1} (t_2\vee _{i_2}(\dots (t_{p-1}\vee _{i_{p-1}} (t_p\vee _{i_p} \vert )))).\end{equation*}
\medskip

\noindent That is

\begin{figure}[h]
\includegraphics[scale=0.6]{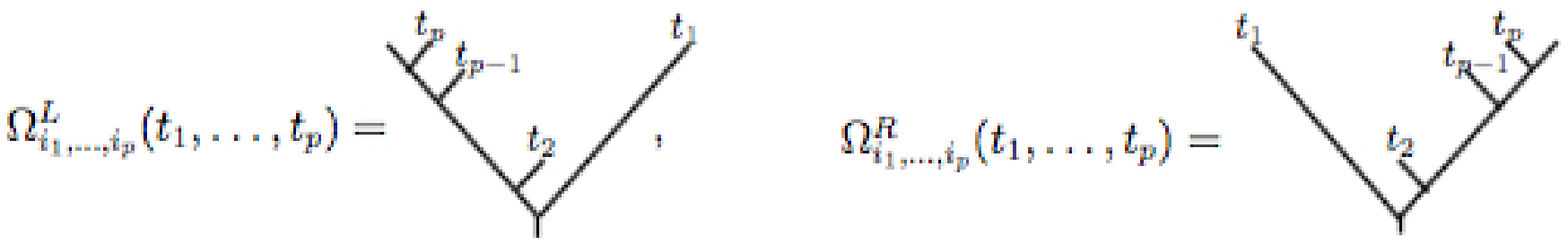}
\end{figure}
\end{enumerate}\end{notation}
\medskip

\begin{remark} \label{coloredcombs} For any tree $t\in {\mathcal Y}_{n-1}^m$ there exist unique non negative integers $p$ and $q$, such that
\begin{equation*} \label{uniquestruct}t = \Omega_{i_1,\dots ,i_p}^L(t_1,\dots ,t_p) =\Omega_{j_1,\dots ,j_q}^R(w_1,\dots ,w_q) ,\end{equation*}
for unique families of colored trees $t_1,\dots, t_p$ and $w_1,\dots ,w_q$, and unique collections of integers $i_1,\dots ,i_p$ and $j_1,\dots ,j_q$ in $\{ 0,\dots ,m\}$, with $i_1=j_1$.

In particular, $t = t^l\vee _{i_1}t^r$, for 
\begin{equation*} t^l = \Omega_{i_2,\dots ,i_p}^L(t_2,\dots ,t_p) = w_1\qquad {\rm and}\qquad t^r= t_1 = \Omega_{j_2,\dots ,j_q}^R(w_2,\dots ,w_q).\end{equation*}\end{remark}

\begin{notation} \label{dimofbasis} Let $b_{m,n}$ denotes the number of elements of the set ${\mathcal B}_n^m$ introduced in Definition \ref{def:basisBm}.
We know that $b_{m,n}$ is the dimension of the subspace of homogeneous elements of degree $n$ of ${\mathfrak D}^m$. Let 
\begin{equation*} f_m(x):=\sum_{n\geq 1} b_{m,n} x^n,\end{equation*}
be the generating series of $\{b_{m,n}\}_{n\geq 1}$, for $m\geq 1$.\end{notation}

\begin{lemma} \label{genseriesofcoloredtrees} The generating series $f_m(x)$ is given by the following conditions\begin{enumerate}
\item $b_{m,1}=1$, for $m\geq 1$,
\item $f_m(x)= x\cdot (1+f_m(x))^{m+1}$.\end{enumerate} \end{lemma}
\medskip

\begin{proo} For $n=1,2$, it is immediate to see that $b_{m,1}=1$ and $b_{m,2}= m$, for $m\geq 1$. 

From Remark  \ref{coloredcombs}, we get that for any colored planar rooted tree $t\in {\mathcal B}_n^m$, there exist unique integers $p\geq 1$, $0\leq i_1<\dots < i_p\leq m$, and unique elements $t_1,\dots ,t_p$ in ${\mathcal B}_m$, such that
\begin{equation*} t= \Omega_{i_1,\dots ,i_p}^L(t_1,\dots ,t_p).\end{equation*}

For $n >2$, the coefficient $b_{m,n}$ of $x^n$ in $f_m(x)$ is 
\begin{align*} \sum_{1\leq p\leq {n-1}}\bigl (\sum_{0\leq i_1<\dots < i_p\leq m} &(\sum_{n_1+\dots +n_p={n-1}}b_{m,n_1}\cdot \ldots \cdot b_{m,n_p})\bigr)=\\
&\sum_{1\leq p\leq {n-1}}\binom{m+1}p (\sum_{n_1+\dots +n_p={n-1}}b_{m,n_1}\cdot \ldots \cdot b_{m,n_p}).\end{align*}

On the other hand, we have that 
\begin{align*} x (1+ f_m(x))^{m+1}=&\\
x \bigl(\sum_{r=0}^{m+1} \binom{m+1}r&\bigl (\sum_{n_1,\dots ,n_r}b_{m,n_1}\cdot \ldots \cdot b_{m,n_p} x^{n_1+\dots +n_r}\bigr)\bigr)=\\
&\sum_{r=0}^{m+1} \binom{m+1}r\bigl (\sum_{n_1,\dots ,n_r}b_{m,n_1}\cdot \ldots \cdot b_{m,n_p} x^{n_1+\dots +n_r+1}\bigr).\end{align*}

So, the coefficient of $x^n$ in $x (1+ f_m(x))^{m+1}$ is
\begin{equation*} \sum_{r=0}^{m+1} \binom{m+1}r \bigl (\sum_{n_1+\dots +n_r+1=n}b_{m,n_1}\cdot \ldots \cdot b_{m,n_p}\bigr ),\end{equation*}
which ends the proof. \end{proo}

\begin{corollary} \label{dimensionofthefree} For $m,n\geq 1$, the dimension of the subspace ${\mathfrak D}_n^m$ of homogeneous elements of degree $n$ of the free $\Dy^m$ algebra ${\mathfrak D}^m$ is equal to $d_{m,n}$.\end{corollary}

The sequence of integers $\{d_{m,n}\}_{n\geq 1}$ describes different families of combinatorial objects, let us mention the planar rooted ${\mbox{}m+1}$-ary trees (see \cite{Sta0} and \cite{HilPed}) and the $m$-Dyck paths (see \cite{Ber}, \cite{BerPre}). 
In Section 4, we describe the free $\Dy^m$ algebra over one element on the space spanned by the set of all $m$-Dyck paths, for $m\geq 1$. The description of the operations $*_i$ on a Dyck are quite technical, but not difficult to explain in terms of the constructions developed in \cite{Ber}, \cite{BerPre} and \cite{BouFusPre}. Clearly, a similar construction may be done in terms of planar rooted $m$-ary trees. However, the last description is less natural, because it requires to cut the original tree
 in three or more subtrees and to glue them again in some complicated way. 

Families of non-symmetric operads $\{{\mathcal P}^m\}_{m\geq 1}$, satisfying that the underlying subspace ${\mathcal P}_n^m$ has dimension $d_{m,n}$, for $n,m\geq 1$, have been defined previously in other authors. Let us mention\begin{enumerate}
\item the operad $^m{\mathcal P}$ introduced by P. Leroux in \cite{Ler}. The operad $^2{\mathcal P}$ coincides with the operad of dendriform algebras and the operad $^3{\mathcal P}$ has dimensions $d_{3,n}$, for $n\geq 1$. For $m\geq 3$, the dimension
 of $^m{\mathcal P}_n$ is different from $d_{m,n}$.
\item the operads of $m$-dendriform algebras, described by J.-C. Novelli in \cite{Nov}
\item the operads ${\mbox{FCat}^{(m)}}$ introduced by S. Giraudo in \cite{Gir}.
\end{enumerate} 

All these operads, as well as $\Dy^m$, are equipped with an associative product, unique up to the product with some scalar.

Let us compare $\Dy^2$ with the operads $^3{\mathcal P}$ (defined in page 5 of \cite{Ler}) and the operad of $2$-dendriform algebras (defined in \cite{Nov}, page 7). Novelli\rq s associative product $*$, as well as the associative product $*$ of 
$^3{\mathcal P}$ must be sent to $\lambda\cdot (*_0+*_1+*_2)$, for some $\lambda \in \K$. So, we may assume that $*\mapsto *_0+*_1+*_2\in \Dy^2$.

Relations $(11)$, $(13)$ and $(16)$ of \cite{Nov}, completely determine the products $\succ $ and $\prec$ on a free $3$-dendriform algebra. Relations $(11)$ and $(16)$ coincide with the relation given in \ref{eq2} for $i=2$ and $i=0$, respectively, 
while relation $(13)$ coincides with the relation \ref{eq1}, for $(i,j)= (0,2)$. So, the product $\succ$ maps to $*_0$, and the product $\prec $ maps to $*_2$.

The product $\circ $ of \cite{Nov} must map to a linear combination of $*_0, *_1$ and $*_2$.  Replacing in relation $(12)$, we get that $\circ \mapsto a *_0 + b *_1$. But, as $*$ maps to $*_0+*_1+*_2$, we get that 
$\circ$ maps to $*_1$. 
In this case, replacing $\succ$, $\circ $ and $\prec$ by $*_0$, $*_1$ and $*_2$, respectively, we see that equations $(14)$ and $(15)$ of \cite{Nov}, state that\begin{enumerate}[(i)]
 \item $(u*_2v)*_1 w = u*_1(v*_1 w + v*_0 w)$,
 \item $(u*_1v+u*_0v) = u*_0 (v*_1 w)$, \end{enumerate}
 which are not valid in a free $\Dy^2$ algebra. 
 
 For $n>2$, the relations $(38)$, $(39)$, $(40)$ and $(42)$ of $m$-dendriform algebras, described in page $15$ of \cite{Nov}, imply that if the products $\succ , \circ_1,\dots ,\circ_{m-1}$ and 
 $\succ$ map to linear combinations of $*_0,\dots ,*_m$, then there exists $0\leq k< m$ such that $\succ \mapsto *_0+\dots +*_k$, $\prec\mapsto *_m$ and $\circ_i$ maps to a linear combination of the products $*_0,\dots ,*_{m-1}$, for $1\leq i <m$. 
A straighforward calculation shows that in this case relation $(41)$ cannot hold in a $\Dy^m$ algebra.
\medskip

Similar arguments hold when we compare the operad $^3{\mathcal P}$ with $\Dy^2$. The associative product $\succ + \prec$ maps to $*_0+*_1+*_2$. As Leroux\rq s products $\succ $ and $\prec$ are uniquely determined on free objects by relations 
$(1), (2)$ and $(3)$, we have two possibilities, to map $\succ$ to $*_0+*_1$ and $\prec$ to $*_2$, or to map $\succ$ to $*_0$ and $\prec$ to $*_1+*_2$. Both choices are equivalent, so we select the first option.

The product $\bullet$ must map to a linear combination of $*_0$, $*_1$ and $*_2$, but using equation $(6)$ we get that $\bullet \mapsto a*_0 + b *_1$, where $a\neq b$. Applying it to equation $(4)$, we obtain 
\begin{equation*}  ((x*_0y)*_0z+ (x*_0y)*_1z) = -a (x*_1 y)*_0z - b (x*_1 y)*_1z,\end{equation*}
which has no solution in a free $\Dy^m$ algebra.
\medskip

Giraudo\rq s operad ${\mbox{FCat}^{(1)}}$ is spanned by products $\bullet_0$ and $\bullet_1$ which satisfy the relations
\begin{equation*} (x\bullet_iy)\bullet_{i+j} z = x\bullet_i(y\bullet_jz),\end{equation*},
for $i,j\geq 0$ and $i+j\leq 1$.

The first equation is $(x\bullet_0y)\bullet_0z = x\bullet_0(y\bullet_0z)$, implies that $\bullet_0$ is associative. So, we we may assume that $\bullet_0\mapsto *_0+*_1$ in the operad $\Dy^1$ of dendriform algebras.

The third equation is $(x\bullet_1y)\bullet_1 z = x\bullet_1(y\bullet_0z)$. If we look for a solution $\bullet_1= a *_0 + b *_1$, we get that $a=0$ and $b=1$. So, the product $\bullet_1$ maps to $*_1$.

Replacing in the second equation $(x\bullet_0y)\bullet_1 z= x\bullet_0(y\bullet _1 z)$, we obtain that $(x*_0 y + x*_1 y)*_1 z = x*_0 (y*_1 z) + x*_1 (y*_1 z)$. This equation implies that $(x *_1 y) *_1 z = x *_1 (y *_1 z)$, which is false in a free 
 $\Dy^1$ algebra. So, the operad ${\mbox{FCat}^{(1)}}$ is not dendriform, and therefore ${\mbox{FCat}^{(1)}}\neq \Dy^1$.  Applying analogous arguments, it is easily seen that ${\mbox{FCat}^{(m)}}\neq \Dy^m$, for $m\geq 1$.
\bigskip

\section {A simplicial object in the category of non-symmetric algebraic operads.} 
\medskip

\subsection{The simplicial structure of $\{\Dy^m\}_{m\geq 0}$}
\medskip

Let ${\mbox{}\Dy ^m{\mbox{-alg}}}$ denotes the category of $\Dy^m$ algebras over $\K$. We want to show that the family of categories $\{{\mbox{}\Dy ^m{\mbox{-alg}}}\}_{m\geq 0}$, where ${\mbox{}\Dy ^0{\mbox{-alg}}}$ is the category of associative algebras, is equipped with face and degeneracy functors, which define a co-simplicial object in the category of small categories.

For any pair of algebraic operads ${\mathcal P}$ and ${\mathcal Q}$, there exists a bijection between functors ${\mathbb F}: {\mbox {}{\mathcal P}{\mbox{-alg}}}\longrightarrow {\mbox {}{\mathcal Q}{\mbox{-alg}}}$ and operad morphisms from ${\mathcal Q}$ to ${\mathcal P}$. So, we get a simplicial object in the category of non-symmetric operads.
\medskip

\begin{definition} \label{defnsimp} Let $(D,*_0,\dots ,*_m)$ be a $\Dy^m$ algebra.\begin{enumerate}
\item Define the product ${\underline {*}}_j$, for $0\leq j\leq m+1$, by
\begin{equation*} {\underline {*}}_j:=\begin{cases} *_j,&{\rm for}\ 0\leq j\leq i{\mbox{}-1},\\
0,&{\rm for}\ j=i,\\
*_{j{-}1},&{\rm for}\ {i+1}\leq j\leq m + 1.\end{cases}\end{equation*}
The binary operations $\{{\underline {*}}_j\}_{0\leq j\leq m + 1}$ define a $\Dy^{m{\mbox{+}1}}$ algebra structure on the vector space $D$. We denote it by ${\mathbb F}_i(D)$. 
\item For $0\leq i\leq m-1$, let ${\overline *}_j$ be the product on $D$ given by
\begin{equation*} {\overline *}_j :=\begin{cases} *_j ,&{\rm for}\ 0\leq j< i,\\
*_i  + *_{i+1},&{\rm for}\ j=i,\\
*_{j-1},&{\rm for}\ i+1<j\leq m,\end{cases}\end{equation*}
The vector space $D$, equipped with the products ${\overline *}_i$, for $0\leq i\leq m-1$, is a $\Dy^{m{\mbox{-}1}}$ algebra, which we denote ${\mathbb S}_i(D)$.\end{enumerate}\end{definition}

\begin{remark} \label{rempropDyckalg} For $m\geq 1$, the map $D\mapsto {\mathbb F}_i(D)$ defines a functor from the category of $\Dy^m$ algebras to the category of $\Dy^{m{\mbox{\it +}1}}$ algebras, for $0\leq i\leq m$. Similarly, the map 
$D\mapsto {\mathbb S}_i(D)$ gives a functor from the category of $\Dy^m$ algebras to the category of $\Dy^{m{\mbox{\it -}1}}$ algebras, for $0\leq i\leq m$. 
 \end{remark}

Remark \ref{rempropDyckalg} implies that ${\mathbb F}_i:\Dy^{m{\mbox{\it +}1}}\longrightarrow \Dy^m$ and 

\noindent ${\mathbb S}_i:\Dy^m\longrightarrow \Dy^{m{\mbox{\it -}1}}$ are operad morphisms, for $0\leq i\leq m$.

The following result is immediate to verify.

\begin{lemma} \label{simpcomp} Let ${\mbox {\it NonSym}}$ be the category of non symmetric algebraic operads. The collection $\{\Dy ^m\}_{m\geq 1}$, equipped with the operad morphisms ${\mathbb F}_i:\Dy^{m{\mbox{\it +}1}}\longrightarrow \Dy^{m}$ and ${\mathbb S}_j:\Dy^{m}\longrightarrow \Dy^{m{\mbox{\it +}1}}$, for $0\leq i\leq m+1$ and $0\leq j\leq m$, is a simplicial object in the category ${\mbox {\it NonSym}}$.\end{lemma}

\subsection{Functors ${\mathbb S}_i$ preserve free objects}
\medskip

For $0\leq i\leq m$, we want to show that the image under the functor ${\mathbb S}_i$ of a free $\Dy ^m$ algebra is free as a $\Dy^{m{\mbox{\it -}1}}$ algebra. From Remark \ref{freeoverone}, it suffices to prove that the image 
${\mathbb S}_i({\mathfrak D}^m)$, of the free $\Dy^m$ algebra over one element, is free as a $\Dy^{m{\mbox{\it -}1}}$ algebra.

In order to do that, we need to introduce a new basis of the vector space ${\mathfrak D}^m$.
\medskip

\begin{definition} \label{def:baseBmk} Given $0\leq k\leq m$, define ${\mathcal B}_n^{m,k}$ to be the set of planar binary rooted trees with $n$ leaves, with the vertices colored by the elements of $\{ 0,\dots ,m\}$ satisfying that, 
for any vertex $v$, the tree $t_v$ fulfills the following conditions\begin{enumerate}[(i)]
\item if $t_v = \Omega_{i_1,\dots ,i_p}^L(t_1,\dots ,t_p) $, with $i_1\neq k$, then either $i_2 = k$ or $i_2 > i_1$,
\item if $t_v = \Omega_{k,\dots ,i_p}^L(t_1,\dots ,t_p) =  \Omega_{k,\dots ,j_q}^R(w_1,\dots ,w_q)$,  then $i_s\in \{ k+1,\dots ,m\}$, for $2\leq s\leq p$, and $j_h\in \{0,\dots ,k\}$, for $2\leq h\leq q$.\end{enumerate}\end{definition}

The basis ${\mathcal B}^m$, introduced in Definition \ref{def:basisBm} coincides with the set ${\mathcal B}^{m,m}$.

Note that the basis ${\mathcal B}^{m,k}$ is obtained from the set of all colored planar binary rooted trees, by replacing certain patterns, in the following way\begin{enumerate}
\item for $0\leq i < j\leq m$ and $i\neq k$, we replace $(t_1\vee_it_2)\vee_jt_3$ by $t_1\vee_i(t_2\vee_jt_3)$,
\item for $k < j\leq m$, we replace $t_1\vee_k(t_2\vee_jt_3)$ by $(t_1\vee_kt_2)\vee_jt_3$,
\item $(t_1\vee_kt_2)\vee_kt_3$ is replaced by $\sum_{j=0}^kt_1\vee_k(t_2\vee_jt_3)- \sum_{i=k+1}^m(t_1\vee_it_2)\vee_kt_3$.\end{enumerate}
\medskip

We first prove that the set ${\mathcal B}^{m,k}$ is a basis of ${\mathfrak D}^m$, for $0\leq k< m$.

\begin{proposition} \label{bijectivemapbasis1} For $n\geq 1$ and $0\leq k\leq m$, there exists a bijective map $\Theta_n: {\mathcal B}_n^m\longrightarrow {\mathcal B}_n^{m,k}$ satisfying that\begin{enumerate}
\item $\Theta_n(t) = t$, for any $t\in {\mathcal B}_n^m\cap {\mathcal B}_n^{m,k}$,
\item when the color of the root of $t$ is $k$, the color of the root of $\Theta_n(t)$ belongs to $\{k,\dots ,m\}$,
\item when the color of the root of $t$ is different from $k$, the color of the root of $\Theta_n(t)$ is the same than the color of the root of $t$,
\item the elements $t$ and $\Theta_n(t)$ describe the same element in ${\mathfrak D}^m$.
\end{enumerate}\end{proposition}
\medskip

\begin{proo} For $n=1$, the map $\Theta_1$ is the identity of ${\mathcal B}_2^m= {\mathcal B}_2^{m,k}$.

For $n=2$, we have that $t\notin {\mathcal B}_2^m\cap {\mathcal B}_2^{m,k}$ if, and only if, $t= \vert\vee_k(\vert\vee_i\vert)$, with $k <i$. In this case, condition \ref{eq1} implies that $t$ represents the same element than 
$(\vert \vee_k\vert)\vee_i\vert$, which belongs to  ${\mathcal B}_2^{m,k}$. So, we define $\Theta (t):= (\vert \vee_k\vert)\vee_i\vert$. 
\medskip

For $n\geq 3$, note that \begin{enumerate}[(i)]
\item if $t=t^l\vee_i t^r \in {\mathcal B}_n^m$ with $i\neq k$, then the color of the root of $\Theta_{n_{\small 1}}(t^l)$ is largest than $i$. So, the element $\Theta_{n_{\small 1}}(t^l)\vee_i\Theta_{n_{\small 2}}(t^r)\in {\mathcal B}_n^{m,k}$. 
Therefore, we get $\Theta_n(t):=\Theta_{n_{\small 1}}(t^l)\vee_i\Theta_{n_{\small 2}}(t^r)$,
\item Suppose that the color of the root of a tree $t\in {\mathcal B}_n^m$ is $k$. If $\Theta_{n_{\small 2}}(t^r)=\Omega_{i_1,\dots ,i_p}^R(w_1,\dots ,w_p)$, with $i_l\leq k$ for $1\leq l\leq p$, then the color of the root of $t^l$ is largest than $k$. Therefore, 
the element $\Theta_{n_{\small 1}}(t^l)\vee_k\Theta_{n_{\small 2}}(t^r)$ belongs to  ${\mathcal B}_n^{m,k}$ and represents the same element than $t$ in ${\mathfrak D}^m$. In this case, we get $\Theta_n(t):= \Theta_{n_{\small 1}}(t^l)\vee_k\Theta_{n_{\small 2}}(t^r)$.\end{enumerate}
\medskip

The unique case which is more complicated is when there exists some $2\leq l\leq p$ such that $k<j_{l}$, where $t=\Omega_{k,j_2,\dots ,j_p}^R(t^l, t_2^r,\dots , t_p^r)$. 

Let $l_0$ be the smallest integer such that $k <j_{\small{l_0}}$. We have that $j_l \leq k$ for $2\leq l< l_0$, and that $t_{j_{\small {l_0}}}^r=\Omega_{s_1,\dots ,s_q}^L(w_1,\dots ,w_q)$, with $j_{l_0}< s_1 <\dots <s_q$ because $t_{j_{\small{l_0}}}^r\in {\mathcal B}^m$.

Using that $k<j_{l_0} < s_1<\dots < s_k$, and applying several times condition \ref{eq1}, we get that $t$ represents the same element that the tree 
\begin{equation*}\Omega_{j_{\small{l_0}},\dots ,j_p}^R(\Omega_{s_1,\dots ,s_q,k,r_1,\dots, r_h}^L(w_1,\dots, w_q, w_{q+1},t_1^l,\dots ,t_h^l), t_{j_{\small{l_0}+1}}^r,\dots ,t_p^r),\end{equation*}
where $w_{q+1}:=\Omega_{j_2,\dots , j_{\small{l_0}-1}}^R(t_2^r,\dots ,t_{j_{\small l_0}-1}^r)$ and $t^l=\Omega_{r_1,\dots ,r_h}^L(t_1^l,\dots ,t_h^l)$.

As $t\in {\mathcal B}_n^m$, we obtain that $k< r_1<\dots <r_h$, and that $j_{l_0} < s_1<\dots < s_q$. Therefore, we get that
\begin{align*} \Theta (t)^l:= \Theta(\Omega&_{s_1,\dots ,s_q,k,r_1,\dots, r_h}^L(w_1,\dots, w_q, w_{q+1},t_1^l,\dots ,t_h^l))= \\
&\Omega_{s_1,\dots ,s_q,k,r_1,\dots, r_h}^L(\Theta(w_1),\dots ,\Theta(w_{q+1}),\Theta(t_1^l),\dots ,\Theta(t_h^l)).\end{align*}

The color of the root of the tree 
\begin{equation*} \Omega_{j_{l_0},j_{l_0+1},\dots ,j_p}^R(\Omega_{s_1,\dots ,s_q,k,r_1,\dots, r_h}^L(w_1,\dots, w_q, w_{q+1},t_1^l,\dots ,t_h^l), t_{j_{\small{l_0}+1}}^r,\dots ,t_p^r)\end{equation*}
is $j_{l_0} > k$, which implies that the tree
$\Omega_{j_{l_0},j_{l_0+1},\dots ,j_p}^R(\Theta (t)^l, \Theta(t_{j_{\small{l_0+1}}}^r),\dots ,\Theta(t_p^r))$ represents the same element than $t$ in ${\mathfrak D}^m$. 
So, we define
\begin{equation*}\Theta (t):= \Omega_{j_{l_0},j_{l_0+1},\dots ,j_p}^R(\Theta (t)^l, \Theta(t_{j_{\small{l_0}+1}}^r),\dots ,\Theta(t_p^r)).\end{equation*}
\medskip

We need to show that $\Theta$ is bijective. For $n=1,2$ the result is immediate. For $n >2$, we apply a recursive argument on $n$. 

Suppose that there exist $\Theta_r^{-1}: {\mathcal B}_r^{m,k}\longrightarrow {\mathcal B}_r^m$, for all $1\leq r<n$, satisfying that\begin{enumerate}[(a)]
\item the composition $\Theta_r\circ \Theta_r^{-1}$ is the identity on ${\mathcal B}_r^{m.k}$,
\item if the root of $t$ is colored by a element smaller or equal than $k$, then the root of $\Theta_r^{-1}(t)$ is colored by the same element.
\item if the root of $t$ is colored by a element $j$ with $k< j$, then the root of $\Theta_r^{-1}(t)$ is colored by an element in $\{ k, j\}$.\end{enumerate}

Let $t=t^l\vee_jt^r\in {\mathcal B}_r^{m,k}$. 

For $j<k$, the root of $t^l$ is colored by $h >j$, thus the root of $\Theta^{-1}(t^l)$ is colored by an element largest than $j$. 

We get that $\Theta^{-1}(t)=\Theta^{-1}(t^l)\vee_j\Theta^{-1}(t^r)$, and it is immediate to see that $\Theta\circ \Theta^{-1}(t)=t$.

For $j=k$, we have that $t^l=\Omega_{s_1,\dots ,s_p}^L(t_1^l,\dots ,t_p^l)$ with $k <s_l$, for $1\leq l\leq p$. The root of $\Theta^{-1}(t^l)$ is colored by an element largest than $k$, which implies that 
\begin{equation*}\Theta^{-1}(t)=\Omega_{s_1,\dots ,s_p}^L(\Theta^{-1}(t_1^l),\dots ,\Theta^{-1}(t_p^l))\vee_k\Theta^{-1}(t^r).\end{equation*} It is easily seen that $\Theta\circ \Theta^{-1}(t)=t$.

For $j>k$, if $t^l=\Omega_{s_1,\dots ,s_p}^L(t_1^l,\dots ,t_p^l)$, for $j <s_1<\dots <s_p$ and $1\leq l\leq p$, then 
\begin{equation*}\Theta_n^{-1}(t) = \Omega_{s_1,\dots ,s_p}^L(\Theta^{-1}(t_1^l),\dots ,\Theta^{-1}(t_p^l))\vee_j\Theta^{-1}(t^r).\end{equation*}

If $s_h=k$, for some $1\leq h\leq p$, then $h$ is unique. As $k<s_{h-1}$, by condition \ref{eq1}, the tree $\Omega_{s_{h-1},k,s_{h+1},\dots ,s_p}^L(t_{h-1}^l,t_h^l,\dots , t_p^l)$ represents the same element than the tree 
\begin{equation*}\Omega_{k,s_{h+1},\dots ,s_p}^L(t_h^l\vee_{s_{h-1}}t_{h-1}^l, t_{h+1}^l,\dots ,t_p^l).\end{equation*}

Since $k<j<s_1<\dots <s_{h-1}$, applying the same argument several times, we get that $t$ represents the same element than the tree
\begin{equation*} \Omega_{k,s_{h+1},\dots ,s_p}^L(w,t_{h+1}^l,\dots ,t_p^l),\end{equation*}
where 
\begin{equation*} w:=\Omega_{j,s_1,\dots ,s_{h-1},r_1,\dots ,r_q}^L(t^r,t_1^l,\dots,t_{h-1}^l,w_1,\dots, w_q),\end{equation*} for $t_h^l=\Omega_{r_1,\dots, r_q}^L(w_1,\dots ,w_q)$.

So, $\Theta_n^{-1}(t):=\Omega_{k,s_{h+1},\dots,s_p}^L(\Theta^{-1}(w),\Theta^{-1}(t_{h+1}^l),\dots ,\Theta^{-1}(t_p^l))$, and the definition of $\Theta$ shows that $\Theta_n\circ\Theta_n^{-1}(t)=t$, which ends the proof.\end{proo}
\medskip

\begin{corollary}\label{cor:basisk} For any $0\leq k\leq m$, the set ${\mathcal B}^{m,k} = {\displaystyle \bigcup _{n\geq 1} {\mathcal B}_n^{m,k}}$ is a basis of the underlying vector space of the free algebra ${\mathfrak D}^m$. \end{corollary}
\medskip

\begin{notation} \label{notn:free} Let $X$ be a set, we denote by ${\mathfrak D}^m(X)$ the free $\Dy^m$ algebra over $X$. The graded set $\bigcup _{n\geq 1}{\mathcal B}_n^{m,k}$ is denoted by ${\mathcal B}^{m,k}$.\end{notation}
\medskip

\begin{lemma} \label{prop:DlDmfree} For any integer $0\leq k < m$, the image of ${\mathfrak D}^m$ under the functor ${\mathbb S}_{k}$ is generated as $\Dy^{m{\mbox{-}1}}$ 
algebra by the graded set ${\mathcal A}^{m,k}$ of all colored trees $t$ in 
${\mathcal B}^{m,k}$, satisfying that either $n = 1$, or $n >1$ and the root of $t$ is colored by $k$. 
\end{lemma}
\medskip

\begin{proo} The $\Dy^{m{\mbox{-}1}}$ algebra structure of ${\mathfrak D}^m$ is given by the products ${\overline {*}}_j =\begin{cases} *_j,& {\rm for}\ 0\leq j< k,\\
*_{k}+*_{k+1},& {\rm for}\ j = k,\\
*_{j-1},& {\rm for}\ k < j < m.\end{cases}$

The underlying vector space of ${\mathbb S}_{k}({\mathfrak D}^m)$ is equal to ${\mathfrak D}^m$. As the set ${\mathcal B}^{m,k}$ is a basis of ${\mathfrak D}^m$ as a $\K$-vector space, it suffices to see that any element in 
${\mathcal B}^{m,k}$ belongs to the $\Dy^{m{\mbox{-}1}}$ algebra generated by the set ${\mathcal A}^{m,k}$, under the operations ${\overline *}_0,\dots , {\overline *}_{m-1}$.

We proceed by induction on the degree $n$. For $n=1,2$, the result is immediate.

For $t = t^l\vee _{i} t^r\in {\mathcal B}_n^{m,k}$, the recursive hypothesis states that the trees $t^l$ and $t^r$ are obtained by applying the products ${\overline {*}}_0,\dots ,{\overline {*}}_{m-1}$ 
to elements of the set ${\mathcal A}^{m,k}$ of degree smaller than $n$. 

We have to analize three different cases\begin{enumerate}
\item for $i < k$, we have that $t = t^l\ {\overline {*}}_i\  t^r$. As $t^l$ and $t^r$ are elements in the $\Dy^{m{\mbox{-}1}}$ algebra generated by ${\mathcal A}^{m,k}$, so is $t$,
\item for $i = k$, as $t\in {\mathcal B}_n^{m,k}$, we get that $t\in {\mathcal A}^{m,k}$,
\item for $i = k+1$, we have that $t = t^l \ {\overline {*}}_{k}\ t^r\ {-}\ t^l *_{k} t^r$ and the root of $t^l$ is colored by $j$, with $j > k+1$ or $j =k$.

As $t^l$ and $t^r$ belong to the free $\Dy^{m{\mbox{\it-}1}}$ algebra ${\mathfrak D}^{m{\mbox{\it-}1}}({\mathcal A}^{m,k})$, the tree $t^l\ {\overline {*}}_{k}\ t^r$ is in ${\mathfrak D}^{m{\mbox{\it-}1}}({\mathcal A}^{m,k})$. 

On the other hand, either $t^l *_{k} t^r\in {\mathcal A}^{m,k}$, or \begin{equation*}t^l *_{k} t^r = (t^l \vee_{k} w^l) \vee _{h} w^r,\end{equation*} for some colored tree $w = w^l \vee_{h} w^r$ and $h > k$.
 
 Applying a recursive argument on the degrees of the elements $t^l \vee_{k} w^l$ and $w^r$ the result follows.
\item For $i > k+1$, we have that $t = t^l \vee_{i-1} t^r$ belongs to ${\mathfrak D}^{m{\mbox{\it-}1}}({\mathcal A}^{m,k})$ by recursive hypothesis.\end{enumerate}
\end{proo}

Lemma \ref{prop:DlDmfree} states that, for any vector space $V$, ${\mathbb S}_{k}({\mathfrak D}^m)$ is a quotient of the free $\Dy^{m{\mbox{\it-}1}}$ algebra ${\mathfrak D}^{m{\mbox{\it-}1}}({\mathcal A}^{m,k})$. For $X$ finite, 
the dimension of the subspace of homogeneous elements of degree $n$ in ${\mathbb S}_{k}(\Dy ^m(X))$ is $d_{m,n}\vert X\vert ^n$. 

So, to prove that ${\mathbb S}_{k}({\mathfrak D}^m(X))$ is isomorphic to ${\mathfrak D}^{m{\mbox{\it-}1}}({\mathcal A}^{m,k}(X))$, it suffices to show that the dimension of the subspace of homogeneous elements of degree $n$ in 
${\mathfrak D}^{m{\mbox{\it-}1}}({\mathcal A}^{m,k})$ is $d_{m,n}$, where ${\mathcal A}^{m,k}$ is the set of trees in ${\mathcal B}_n^m$ with the vertices colored by $\{0,\dots , m\}$ and the root colored by $k$.
\medskip

Recall that, for any graded vector space $V=\bigoplus_{n\geq 1}V_n$ such that each $V_n$ is finite dimensional, the generating series of $V$ is $v(x) := {\displaystyle \sum_{n\geq 1} {\mbox{dim}_{\K}(V_n)} x^n}$.

\begin{lemma}\label{lemform}
Let $d_m(x)$ be the generating series of the free $\Dy^m$ algebra ${\mathfrak D}^m$. We have that
\begin{equation} \label{eq5} d_m(x)=d_k(x\cdot (1+d_m(x)^{m-k}),\end{equation}
for all $0\leq k\leq m$.
\end{lemma}
\medskip

\begin{proo} Clearly, it suffices to prove the result for $k=m-1$. Let $g_m(x)$ be the inverse series of $d_m(x)$ ($g$ exists because $d(0)=0$). 

Since $x\cdot (1+d_m(x))^{m+1}=d_m(x)$, replacing $x$ by $g_m(x)$ we obtain that
\begin{equation*}g_m(x)=\displaystyle\frac{x}{(1+x)^{m+1}},\end{equation*}
which implies that  $(1+x)\cdot g_m(x) = g_{m-1}(x)$. So, replacing $x$ by $d_m(x)$ and applying $d_{m-1}(x)$ to both sides, we get the desired formula
\begin{equation*}d_{m-1} (x\cdot (1+d_m(x))) = d_m(x).\end{equation*}
\end{proo}

Applying Lemmas \ref{prop:DlDmfree} and \ref{lemform}, we get the following result.

\begin{proposition} \label{pro:freenessform} The $\Dy^{m{\mbox{\it-}1}}$ algebra ${\mathbb S}_{k}({\mathfrak D}^m(X))$ is free, for any $0\leq k\leq m-1$.\end{proposition}
\medskip

\begin{proo} Applying Lemmas \ref{prop:DlDmfree} and \ref{lemform}, it suffices to prove that the number of elements in ${\mathcal A}_n^{m,k}$ is $d_{m,n-1}$, for $0\leq k\leq m$. 

The number of elements of ${\mathcal B}_{n-1}^{m,k}$ is $d_{m,n-1}$. To end the proof we define a bijective map $\theta _n$ from ${\mathcal B}_{n-1}^{m,k}$ to ${\mathcal A}_n^{m,k}$, for $n\geq 2$. 

For $n = 2$, $\theta _1(\vert )$ is the unique planar binary rooted tree with two leaves and the root colored by $k$.

Let $t =t^l\vee_{h}t^r$ be an element of ${\mathcal B}_{n-1}^{m,k}$. \begin{enumerate}
\item For $h > k$, let $t = \Omega_{h,i_2,\dots,i_p}^L(t^r,t_2,\dots ,t_p)$.\begin{enumerate}
\item If $i_p >\dots > i_2 > h > k$, then 
$\theta _n(t) := t\vee _{k} \vert .$
\item If there exists one integer $1\leq s\leq p$ such that $i_s = k$, then the $s$ is unique and $\theta _n(t)$ is defined by the formula
\begin{equation*}\theta _n(t) :=\Omega_{h,i_2,\dots, i_{s{-}1}}^L(t^r,t_2,\dots ,t_{s{-}1})\vee_{k}\Omega_{k,i_{s+1},\dots ,i_p}^L(t_s,\dots, t_p).\end{equation*}\end{enumerate}
\item For $h\leq k$, let $t = \Omega _{h,j_2,\dots ,j_q}^R(t^l,w_2,\dots ,w_q)$.\begin{enumerate}
\item If $j_i\leq k$ for any $2\leq i\leq q$, then 
$\theta _n(t) := \vert \vee _{k} t$.
\item Otherwise, let $2\leq s\leq q$ be the minimal integer satisfying that $j_s >k$. As $t\in {\mathcal B}^{m,k}$, we know that $k\notin \{h,j_1,\dots ,j_{s-1}\}$.
Define $\theta _n(t)$ as the element
\begin{equation*}\theta _n(t):= \Omega_{j_s,\dots ,j_q}^R(w_s,\dots ,w_q)\vee_{k}\Omega_{h,j_2,\dots ,j_{s{-}1}}^R(t^l,w_2,\dots ,w_{s{-}1}).\end{equation*}
\end{enumerate}
\end{enumerate}

It is not difficult to verify that $\theta _n$ is bijective for all $n\geq 2$. So, the result is proved.
\end{proo}

The following result is a straightforward consequence of Proposition \ref{pro:freenessform}. 

\begin{theorem} The image of a free $\Dy ^m$ algebra ${\mathfrak D}^m(X)$ under the functor ${\mathbb S}_{k}$ is a free $\Dy^{m{\mbox{-}1}}$ algebra, for $0\leq k<m$.\end{theorem}
\medskip

Note that, by composing the functors ${\mathbb S}_{k}$, we get that the associative algebra $({\mathfrak D}^m(X), *_0+\dots +*_m)$ is free, for any set $X$.
\bigskip

\section{Dendriform posets} 
\medskip 

Let ${\mathcal Y}:=\{{\mathcal Y}_n,\leq_{Ta}\}_{n\geq 1}$ be the family of sets ${\mathcal Y}_n$ of planar binary rooted trees. The vector space $V_{\small{\mathcal Y}}:=\bigoplus_{n\geq 1}\K[{\mathcal Y}_n]$
admits a structure of dendriform algebra, satisfying that the products $\succ $, $\prec$ and $* = \succ +\prec$ are defined in terms of intervals of the Tamari order. This result was described in \cite{LodRon1}. In the same work, an analogous result 
was proved for the family ${\bold \Sigma}=\{\Sigma_n\}_{n\geq 1}$, where $\Sigma_n$ is the set of permutations of $n$ elements, where the Tamari order is replaced by the weak Bruhat order (see also \cite{Foi}). Both results were extended in \cite{PalRon} to the families of all surjective maps and of all planar rooted trees, where the partial orders of these sets are generalizations of the weak Bruhat order and of the Tamari order, respectively. 

We describe the conditions that a family of partially ordered sets ${\bold P}=\{P_n\}_{n\geq 1}$ must fulfill in order to\begin{enumerate}
\item get a natural structure of dendriform algebra on the vector space $V_{\bold P}=\bigoplus_{n\geq 1}\K[P_n]$, spanned by the graded set $\bigcup _{n\geq 1}P_n$.
\item get an $\Dy^m$ algebra structure on the vector space spanned by the set $\bigcup_{n\geq 1}{\mbox {Simp}(P_n)^m}$ of $m$-simpleces of the partially ordered sets $P_n$, $n\geq 1$. \end{enumerate}

This type of families of partially ordered sets provide examples of $\Dy^m$ algebras. 

\subsection{Basic constructions} 
\medskip

\begin{definition} \label{defdendposet} A {\it dendriform poset} is a family of partially orederd sets ${\bold P}=\{P_n\}_{n\geq 1}$, equipped with four binary graded products $/$, $\perp$, $\top$ and $\backslash$, satisfying the following conditions\begin{enumerate}
\item for $n,r\geq 1$, the maps $/, \perp, \top, \backslash : P_n\times P_r\longrightarrow P_{n+m}$, preserve the orders (where the partial order of $P_n\times P_m$ is the componentwise one),
\item for any pair of elements $x\in P_n$ and $y\in P_m$, the interval $[ x/ y;x\backslash y]$ is the disjoint union of the intervals $[ x/y;x\perp y]$ and $[x\top y;x\backslash y]$,
\item for any elements $x\in P_n$, $y\in P_r$ and $z\in P_s$, there exists a bijective order preserving map $\varphi_{x,y,z}$ from the set 
\begin{equation*} L(x,y,z) :=\{ (u,v)\in P_{n+r}\times P_s\mid x/y\leq u\leq x\backslash y\ ,\ u/z\leq v\leq u\backslash z\},\end{equation*}
to the set 
\begin{equation*} R(x,y,z) :=\{ (u,v)\in P_{n}\times P_{r+s}\mid y/z\leq u\leq y\backslash z\ ,\ x/u\leq v\leq x\backslash u\},\end{equation*}
satisfying that\begin{enumerate}
\item the restriction of $\varphi_{x,y,z}$ to the set 
\begin{equation*} L^{\succ}(x,y,z):=\{ (u,v)\in P_{n+r}\times P_s\mid x/y\leq u\leq x\backslash y\ ,\ u/z\leq v\leq u\perp z\},\end{equation*}
 gives a bijection with the set 
 \begin{equation*} R^{\succ}(x,y,z):=\{ (u,v)\in P_{n}\times P_{r+s}\mid y/z\leq u\leq y\perp z\ ,\ x/u\leq v\leq x\perp u\},\end{equation*}
\item the restriction of $\varphi_{x,y,z}$ to the set 
\begin{equation*}L^{\prec}(x,y,z):=\{ (u,v)\in P_{n+r}\times P_s\mid x\top y\leq u\leq x\backslash y\ ,\ u\top z\leq v\leq u\backslash z\}\end{equation*} gives a bijection with the set 
\begin{equation*}R^{\prec}(x,y,z):=\{ (u,v)\in P_{n}\times P_{r+s}\mid y/z\leq u\leq y\backslash z\ ,\ x\top u\leq v\leq x\backslash u\}\end{equation*},\end{enumerate}
\item if there exist  $u\leq  v$ in $P_{n+r}$ satisfying that $x/y\leq u\leq x\backslash y$ and $z/w\leq v\leq z\backslash w$, for some elements $x, z\in P_n$ and $y,w\in P_r$, then $x\leq z$ in $P_n$ and $v\leq w$ in $P_r$.
\item if there exist $u\leq  v$ in $P_{n+r}$ satisfying that $x/y\leq u\leq x\perp y$ and $x\top y\leq v\leq x\backslash y$, for some elements $x\in P_n$ and $y\in P_r$, then $v\not\leq u$.\end{enumerate}\end{definition}
\medskip

\begin{notation} Let ${\bold P}=\{P_n\}_{n\geq 1}$ be a dendriform poset. We denote by $V_{\small{\bold P}}:=\bigoplus_{n\geq 1}\K[P_n]$ the graded vector space spanned by ${\bold P}$. \end{notation}
\medskip

\begin{definition} \label{dendrstructdendposet} Let ${\bold P}=\{P_n\}_{n\geq 1}$ be a dendriform poset. The products $\succ $ and $\prec $ are defined by
\begin{align} x\succ y &:=\sum_{x/y\leq u\leq x\perp y} u,\\
x\prec y&:=\sum_{x\top y\leq u\leq x\backslash y} u,\end{align}
for any elements $x, y \in \bigcup_{n\geq 1}P_n$. Both products are extended by linearity to the vector space $V_{\small{\bold P}}$.\end{definition}
\medskip

\begin{proposition} \label{dendestposet} Let ${\bold P}=\{P_n\}_{n\geq 1}$ be a dendriform poset. The vector space $V_{\small{\bold P}}$ equipped with the products $\succ $ and $\prec$ is a dendriform (or $\Dy^1$) algebra. Conversely, if  
$(V_{\small{\bold P}}, \succ ,\prec)$ is a dendriform algebra, then ${\bold P}$ satisfies conditions $(2)$ and $(3)$ of Definition \ref{defdendposet}.\end{proposition}
\medskip

\begin{proo} Condition $(2)$ is equivalent to $x*y=x\succ y + x\prec y =\sum_{x/y\leq u\leq x\backslash y} u$, for any elements $x, y\in V_{\small{\bold P}}$.

Condition $(3)$ is equivalent to\begin{enumerate}
\item $x*(y*z)= (x*y)*z$.
\item $x\succ (y\succ z) = (x*y)\succ z$,
\item $x\prec (y*z) = (x\prec y)\prec z$,\end{enumerate}
for any elements $x,y,z\in \bigcup_{n\geq 1}P_n$, therefore for any $x, y$ and $z$ in $V_{\small{\bold P}}$.

As $*=\succ + \prec$ is associative, we get that 
\begin{align*} x\succ (y\prec z) &= x*(y*z) - x\prec (y*z)-x\succ (y\succ z)=\\
&(x*y)*z - (x\prec y)\prec z -(x*y)\succ z =(x\succ y)\prec z,\end{align*}
for any elements $x,y$ and $z$ in $V_{\small{\bold P}}$. 

So, $(V_{\small{\bold P}}, \succ ,\prec)$ is a dendriform algebra.\end{proo} 
\bigskip

\subsection{Examples} {\bf 1) The facial order on surjective maps} 
\begin{notation}\label{notperms} Let ${\mbox{Surj}_n^r}$ denotes the set of surjective maps from $[n]$ to $[r]$, for $1\leq r\leq n$. For any $f\in {\mbox{Surj}_n^r}$, we denote it by the tuple of its images $f=(f(1),\dots ,f(n))$. Note that ${\mbox{Surj}_n^n}$ is the group $\Sigma_n$ of permutations of $n$ elements, for $n\geq 1$. 

We denote the disjoint union $\bigcup_{r=1}^n{\mbox{Surj}_n^r}$ by ${\mbox{Surj}_n}$.

For $1\leq i\leq {n-1}$, let $s_i$ be the element of $\Sigma _n$ which exchanges $i$ and $i+1$. The surjective map $\tau_i\in {\mbox{Surj}_n^{n{\mbox{-}1}}}$ is the map given by 
\begin{equation*}\tau_i(j):=\begin{cases}j,&{\rm for}\ 1\leq j\leq i,\\
j-1,&{\rm for}\ 1\leq j\leq i,\end{cases}.\end{equation*}

Let $S$ and $T$ be two disjoint subsets of $[n]$, we say that $S <T$ if $s < t$, for any pair of elements $s\in S$ and $t\in T$.

For any $S\subseteq [n]$ and any $f\in {\mbox{Surj}_n}$, let $f\vert_S$ denote the restriction of $f$ to the set $S$. Suppose that $S=\{j_1,\dots ,j_s\}$ and that $r$ is a positive integer, we denote by $S+r$ the set $\{j_1+r,\dots ,j_s+r\}$.\end{notation}

\begin{definition} \label{standadr} Let $f:[n]\longrightarrow {\mathbb N}$ be a map, the {\it standardization} of $f$ is the unique surjective map ${\mbox{std}(f)}:[n]\longrightarrow [\vert {\mbox {Im}(f)}\vert ]$ satisfying that ${\mbox{std}(f)}(i) < {\mbox{std}(f)}(j)$ if, and only if, $f(i) < f(j)$, for all elements $1\leq i,j\leq n$.\end{definition}

\begin{definition} \label{defn:facialorder} The facial order on the set ${\mbox{Surj}_n}$ is the transitive relation $\leq_{fa}$ spanned by the following covering relations \begin{enumerate}
\item If $f^{-1}(i) < f^{-1}(i+1)$, then $f \lessdot \tau_i\circ f$, for $1\leq i <\vert{\mbox {Im}(f)}\vert $,
\item if $f^{-1}(i)=\{j_1<\dots < j_s\}$, then $f\lessdot f_k$, where $f_k$ is the map 
\begin{equation*} f_k(j):=\begin{cases}f(j)&{\rm for}\ 1\leq f(j)<i,\\
i&{\rm for}\ f(j)=i\ {\rm and}\ j\in\{j_{k+1},\dots ,j_s\}\\
i+1&{\rm for}\ f(j)=i\ {\rm and}\ j\in \{j_1,\dots ,j_k\}\\
f(j)+1&{\rm for}\ f(j)>i,\end{cases}\end{equation*}
 for any $1\leq k<s$.\end{enumerate}\end{definition}
\medskip

In order to define a dendriform poset structure on ${\bold{\mbox{Surj}}}$, we must define the graded products $/$, $\perp$, $\top$ and $\backslash$, these operations were introduced in \cite{PalRon}.

\begin{definition} \label{operforperm} Let $f\in {\mbox{Surj}_n^s}$ and $g\in {\mbox{Surj}_r^h}$ be two surjective maps, define
\begin{align*} f/g&:=(f(1),\dots ,f(n),g(1)+s,\dots ,g(r)+s),\\
f\perp g(j)&:=\begin{cases}f(j)+{\mbox{}h-1},&{\rm for}\ 1\leq j\leq n,\\
g(j-n),&{\rm for}\ n+1\leq j\leq n+r\ {\rm and}\ g(j)<h,\\
s+h,&{\rm for}\ g(j)=h.\end{cases}\\
f\top g(j)&:=\begin{cases}f(j),&{\rm for}\ f(j)<s,\\
s+h,&{\rm for}\ f(j)=s,\\
g(j-n)+h,&{\rm for}\ n<j\leq n+r,\end{cases}\\
f\backslash g(j)&:=(f(1)+h,\dots, f(n)+h,g(1),\dots ,g(r)).\end{align*}\end{definition}
\medskip

\begin{proposition} \label{propsurjdendposet} The family ${\mbox{\bf Surj}}=(\{{\mbox{Surj}_n}\}_{n\geq 1}, \leq_{fa})$ is a dendriform poset.\end{proposition}
\medskip

\begin{proo} A recursive argument on $n\geq 1$ shows that the products $/, \perp ,\top$ and $\backslash : {\mbox{Surj}_n}\times {\mbox{Surj}_r}\longrightarrow {\mbox{Surj}_{n+r}}$ preserve the order.

The proof that $V_{\small{\mbox {\bf Surj}}}$ is a dendriform algebra with the products 
$f\succ g:=\sum _{f/g\leq_{fa} u\leq_{fa}f\perp g} u$ and $f\prec g:=\sum _{f\top g\leq_{fa} u\leq_{fa}f\backslash g} u$, and that the associative product $*:=\succ +\prec$ is given by $f*g:=\sum _{f/g\leq_{fa} u\leq_{fa}f\backslash g} u$ was done in \cite{PalRon}, we refer to this reference for the details of the proof. This result implies that ${\bold{\mbox{Surj}}}$ fulfills the second and third conditions of Definition \ref{defdendposet}. We need to prove that it also satisfies the last two ones.
\medskip

Note that, for any pair of elements $f\in {\mbox{Surj}_n}$ and $g\in {\mbox{Surj}_r}$, we have that ${\mbox{std}((f\odot g)\vert_{\{1,\dots ,n\}}}=f$ and ${\mbox{std}((f\odot g)\vert_{\{n+1,\dots ,n+r\}}}=g$, where $\odot$ is any element of the set $\{/,\perp , \top , \backslash\}$.\medskip

On the other hand, for $f\leq _{fa} g$ in ${\mbox{Surj}_n^s}$ and any subset $S\subseteq [r]$, it is immediate to see that $f\vert_S\leq_{fa}g\vert_S$. 

If $u\leq_{fa}w$ in ${\mbox{Surj}_{n+r}}$ and there exist $f_1,f_2\in {\mbox{Surj}_n}$ and $g_1,g_2\in {\mbox{Surj}_r}$ such that
\begin{equation*} f_1/g_1\leq_{fa} u\leq _{fa}f_1\backslash g_1\ {\rm and}\ f_2/g_2\leq_{fa}w\leq_{fa} f_2\backslash g_2,\end{equation*}
then $f_1={\mbox{std}(u\vert_{\{1,\dots ,n\}})}\leq_{fa} {\mbox{std}(w\vert_{\{1,\dots ,n\}})}=f_2$ and $g_1={\mbox{std}(u\vert_{\{n+1,\dots ,n+r\}})}\leq_{fa} {\mbox{std}(w\vert_{\{n+1,\dots ,n+r\}})}=g_2$. Therefore, $f_1\leq_{fa} f_2$ and $g_1\leq_{fa}g_2$, and the fourth condition is satisfied.
\medskip

Let $f_1,f_2\in {\mbox{Surj}_n}$ and $g_1,g_2\in {\mbox{Surj}_r}$, and let $u$ and $w$ be elements in ${\mbox{Surj}_{n+r}}$ satisfying that
\begin{equation*} f_1/g_1\leq_{fa} u\leq_{fa}f_1\perp g_1\ {\rm and}\ f_2\top g_2\leq_{fa} w\leq_{fa}f_2\backslash g_2.\end{equation*}
We need to prove that $w\not \leq_{fa} u$. 

Suppose that $u\in {\mbox{Surj}_{n+r}^k}$ and $w\in {\mbox{Surj}_{n+r}^h}$ satisfy that $w\leq_{fa}u$. The definition of $\leq_{fa}$ implies that, for any $i\in w^{-1}(k)$ there exists $1\leq j_i\leq i$ such that $u(j_i)=k$.

On the other hand, suppose that $f_i\in{\mbox{Surj}_n^{s_i}}$ and $g_i\in {\mbox{Surj}_r^{h_i}}$, for $i=1,2$. We have that\begin{enumerate}[(a)]
\item $f_1/g_1$ and $f_1\perp g_1$ belong to ${\mbox{Surj}_{n+r}^{s_1+h_1}}$ satisfy 
\begin{equation*} (f_1/g_1)^{-1}(s_1+h_1)=(f_1\perp g_1)^{-1}(s_1+h_1)= g_1^{-1}(h_1)+n,\end{equation*}
which implies that $u^{-1}(k) = g_1^{-1}(h_1)+n$.
\item $f_2\top g_2$ and $f_2\backslash g_2$  belong to ${\mbox{Surj}_{n+r}^{s_2+h_2}}$, and satisfy\begin{enumerate}[(i)] 
\item $(f_2\top g_2)^{-1}(s_2+h_2)= f_2^{-1}(s_2)\bigcup (g_2^{-1}(h_2)+n),$
\item $(f_2\backslash g_2)^{-1}(s_2+h_2)= f_2^{-1}(s_2),$\end{enumerate}
which imply that
\begin{equation*}f_2^{-1}(s_2)\subseteq w^{-1}(l) \subseteq   f_2^{-1}(s_2)\bigcup (g_2^{-1}(h_2)+n).\end{equation*}\end{enumerate}

Point $(a)$ implies that $u^{-1}(k)\subseteq \{n+1,\dots , n+r\}$. Let $i\in f_2^{-1}(s_2)$, by $(b)$ we know that $w(i)=l$. So, there exists $1\leq j_i\leq i\leq n$ such that $u(j_i)=k$, which is false. So, $w\not\leq_{fa}u$, and the last condition of Definition \ref{defdendposet} is fulfilled.\end{proo}
\medskip

\begin{remark} \label{facialBruhat} The symmetric group $\Sigma_n$ is generated by the permutations $s_i$, for $1\leq i\leq {\mbox{}n-1}$. The length of a permutation $\sigma\in \Sigma_n$ is the minimal integer $l(\sigma)$ satisfying that $\sigma=s_{j_1}\cdot \dots \cdot s_{j_{l(\sigma)}}$. The left weak Bruhat  order $\leq_B$ on $\Sigma_n$ is the transitive relation spanned by
\begin{equation*} \sigma \lessdot s_i\cdot \sigma, \ {\rm whenever}\ l(s_i\cdot \sigma) =l(\sigma )+1.\end{equation*}
The restriction of the order $\leq_{fa}$ to the subset $\Sigma_n={\mbox{Surj}_n^n}$ gives the left weak Bruhat order.\end{remark}

Proposition \ref{propsurjdendposet} and Remark \ref{facialBruhat} imply that the family ${\bold \Sigma}=\{\Sigma_n\}_{n\geq 1}$ is also a dendriform poset, via the inclusion map $\Sigma_n\subseteq {\mbox{Surj}_n}$. It suffices to note that ${\bold \Sigma}$ is closed under the products $/$. $\perp$ and $\backslash$, and that the product $\top$ is given by
\begin{equation*} \sigma \top \tau (j) = \begin{cases} f(j),&{\rm for}\ f(j)<n,\\
n+r,&{\rm for}\ f(j)=n,\\
g(j-n)+{\mbox{}n-1}, &{\rm for}\ n<j<n+r.\end{cases}\end{equation*}
\bigskip

{\bf 2) The partial order on planar rooted trees} 
\medskip

Let ${\mathcal T}_n^r$ be the set of planar rooted trees with $n+1$ leaves and $n-r$ internal vertices, for $0\leq r\leq n-1$. Note that ${\mathcal T}_n^0={\mathcal Y}_n$ is the set of planar binary rooted trees with $n+1$ leaves.

\begin{notation} \label{planartrees} We denote by ${\mathcal T}_n$ the set ${\displaystyle \bigcup _{r=0}^{n-1}{\mathcal T}_n^r}$. For a collection of planar rooted trees $t^1,\dots ,t^p$, $\bigvee (t^1,\dots ,t^p)$ denotes the planar tree obtained by joining the roots of the trees $t^1,\dots , t^p$, disposed from left to right, to a new root. If $t^i\in {\mathcal T}_{n_i}^{r_i}$, for $1\leq i\leq p$, then $\bigvee (t^1,\dots ,t^r)\in {\mathcal T}_{n}^{r}$, where $n+1=n_1+\dots +n_p+p$ and $r+2=r_1+\dots+r_p+p$.

For $t\in {\mathcal T}_n$, we denote $\vert t\vert = n$.\end{notation} 

There exist surjective maps (see \cite{NovThi0}, \cite{PalRon}) $\Gamma_n: {\mbox{Surj}_n}\longrightarrow  {\mathcal T}_n$, satisfying that the inverse image $\Gamma_n^{-1}(t)$ is an interval in the facial order 
of ${\mbox{Surj}_n}$, for $n\geq 1$. The facial order $\leq_{fa}$ of the set ${\mbox{Surj}_n}$ induces a partial order $\leq _T$ on the set ${\mathcal T}_n$, for $n\geq 1$. 

These partially ordered sets were studied in \cite{PalRon}, in order to show that the dendriform (and tridendriform) algebra structure, defined in \cite{LodRon2} on the graded vector space 
$\K[{\mathcal T}_{\infty}]:=\bigoplus_{n\geq 1} \K[{\mathcal T}_n]$, may be described in terms of intervals of these orders. We briefly describe the main results needed to show that the family $(\{{\mathcal T}_n\}_{n\geq 1},\leq_T)$ is a dendriform poset, for a more detailed description we refer to \cite{PalRon}.

The following result was proved in \cite{PalRon}. 

\begin{lemma} \label{partordertrees} The facial order of ${\mbox{Surj}_n}$ induces a partial order $\leq_T$ on ${\mathcal T}_n$, which extends the Tamari order defined on the set ${\mathcal Y}_n$ of planar binary rooted trees. The partial order $\leq_T$ is completely determined by the following conditions\begin{enumerate}
\item for $n=2$, the order is described by 
\begin{figure}[h]
\includegraphics[scale=0.4]{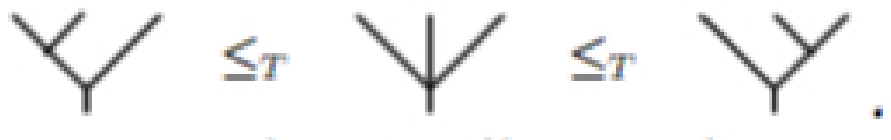}
\end{figure}
\item for a family of trees $t^1,\dots ,t^p$, suppose that $t^i=\bigvee(t^{i1},\dots ,t^{is})$, with $s>1$, for some $1\leq l<p$. In this case
\begin{equation*} \bigvee(t^1,\dots ,t^p)<_T \bigvee (t^1,\dots ,t^{i-1},t^{i1}, \dots ,t^{is},t^{i+1},\dots ,t^p).\end{equation*}
\item for a family of trees $t^1,\dots ,t^p$, with $p\geq 3$, 
\begin{equation*} \bigvee(t^1,\dots ,t^p)<_T \bigvee (t^1,\dots ,t^{i-1},\bigvee(t^i,\dots, t^{j}),t^{j+1},\dots ,t^p),\end{equation*}
for any pair $0< i < j < p$.
\item given two families of trees $t^1,\dots ,t^p$ and $w^1,\dots, w^p$, with $p\geq 2$, satisfying that $t^l\leq_T w^l$, for $1\leq l\leq p$, 
\begin{equation*}\bigvee(t^1,\dots ,t^p)\leq_T\bigvee(w^1,\dots ,w^p).\end{equation*}\end{enumerate}\end{lemma}
\medskip

\begin{definition} \label{slicing} For any tree $t\in {\mathcal T}_n$ and any $0\leq l\leq n$, the {\it $i^{th}$ restriction} of $t$ is the pair $(t_{(1)}^l,t_{(2)}^l)\in {\mathcal T}_l\times {\mathcal T}_{n-l}$ recursively defined by\begin{enumerate}[(a)]
\item for $i=0$, $(t_{(1)}^0,t_{(2)}^0):=(\vert ,t)$,
\item for $i=n$, $(t_{(1)}^n,t_{(2)}^n):=(t, \vert )$,
\item for $n=2$ and $i=1$, $(t_{(1)}^1,t_{(2)}^1):= (\bigvee(\vert ,\vert), \bigvee(\vert ,\vert))$, for any $t\in {\mathcal T}_2$,
\item For $t=\bigvee (t^1,\dots ,t^p)$, with $\vert t^j\vert =n_j$, and $1< l+1< n_1+\dots +n_p+p$, define
\begin{equation*} ((t_{(1)}^l,t_{(2)}^l):= (\bigvee (t^1,\dots ,(t^{k})_{(1)}^{l-(n_1+\dots +n_k+k)}), \bigvee ((t^{k})_{(2)}^{l-(n_1+\dots +n_k+k)}, t^{k+1},\dots ,t^p)),\end{equation*}
where $1\leq k\leq p$ is the unique integer such that $n_1+\dots +n_{k-1}+k\leq l+1< n_1+\dots +n_k+k+1$.\end{enumerate}\end{definition}
\medskip

In fact the $i^{th}$ restriction of a tree $t$ is the result of slicing in two the tree $t$ following the path from the $i^{th}$ leaf to the root. Moreover, a long but simple recursive argument shows the following result.

\begin{lemma} \label{lemorder} For any tree pair of trees $t\leq_T w$ in ${\mathcal T}_n$ and any $0\leq j\leq n$, we have that $t_{(i)}^j\leq_Tw_{(i)}^j$, for $i=1,2$.\end{lemma}

\begin{definition} \label{opontrees} Let $t\in {\mathcal T}_n$ and $w\in {\mathcal T}_r$ be two planar rooted trees. The planar rooted trees $t/w$ and $t\backslash w$ in ${\mathcal T}_{n+r}$ are defined as follows\begin{enumerate}
\item $t/w$ is the tree obtained by grafting the root of $t$ to the leftmost leaf of $w$,
\item $t\backslash w$ is the tree obtained by grafting the root of $w$ to the rightmost leaf of $t$.\end{enumerate}\end{definition}
\medskip

\begin{remark} \label{slicing} For $t\in {\mathcal T}_n$ and $w\in {\mathcal T}_r$, we have that $((t/w)_{(1)}^n, (t/w)_{(2)}^n) = (t,w)=((t\backslash w)_{(1)}^n, (t\backslash w_{(2)}^n))$.\end{remark}
\medskip

\begin{proposition} \label{propsurjdendposet} The family ${\bold{\mathcal T}}:=(\{{\mathcal T}_n\}_{n\geq 1}, \leq_{T})$ is a dendriform poset.\end{proposition}
\medskip

\begin{proo} In \cite{LodRon2} three binary operations $\succ$, $\cdot$ and $\prec$ are defined on the vector space $\K[{\mathcal T}_{\infty}]$, spanned by the set of all planar rooted trees, in such a way that the products $*_0:=\succ$ and 
$*_1:= \cdot + \prec$ define a dendriform (or $\Dy^1$) algebra. 

For any pair of planar rooted trees $t=\bigvee (t^1,\dots ,t^p)$ and $w=\bigvee (w^1,\dots ,w^q)$, define the trees
\begin{align*} t\perp w&:=\bigvee (t\backslash w^1,w^2,\dots ,w^q)\\
t\top w&:=\bigvee (t^1,\dots ,t^p/w^1,\dots ,w^q),\end{align*}
where $/$ and $\backslash$ are the products described in Definition \ref{opontrees}.
The maps $/, \perp ,\top$ and $\backslash$ are morphisms of partially ordered sets from ${\mathcal T}_n\times {\mathcal T}_r$ to ${\mathcal T}_{n+r}$, for $n,r\geq 1$. Moreover, the products $*_0$ and $*_1$ are given by the equalities
\begin{equation*} t *_0 w = \sum_{t/w\leq_T u\leq_T t\perp w}u,\qquad {\rm and}\qquad
t *_1w = \sum_{t\top w\leq_T u\leq_T t\backslash w}u,\end{equation*}
 for any planar rooted trees $t$ and $w$. We refer to \cite{PalRon} for the details of the proof. 
 
The previous result implies that the family ${\bold{\mathcal T}}$, satisfies the first three conditions of Definition \ref{defdendposet}.
\medskip

Suppose that $t_1/w_1\leq_T u\leq_T t_1\backslash w_1$ and $t_2/w_2\leq_T v\leq_T t_2\backslash w_2$ are such that $u\leq_T v$, for elements $t_1, t_2\in {\mathcal T}_n$ and $w_1, w_2\in {\mathcal T}_r$. 
Applying Lemma \ref{lemorder} and Remark \ref{slicing}, we get that 
\begin{equation*} t_1=u_{(1)}^n\leq_T v_{(1)}^n=t_2.\qquad {\rm and}\qquad w_1=u_{(2)}^n\leq_T w_{(2)}^n = w_2,\end{equation*} which implies that condition $(4)$ of Definition \ref{defdendposet} is satisfied.
\medskip

To end the proof we need to see that if $t_1/w_1\leq_T u\leq_T t_1\perp w_1$ and $t_2\top w_2\leq _T v\leq _T t_2\backslash w_2$, for two pairs of elements $t_1, t_2\in  {\mathcal T}_n$ and $w_1, w_2\in {\mathcal T}_r$, then $v\not\leq_T u$. 

Suppose that $t_i=\bigvee(t_i^1,\dots ,t_i^{p_i})$ and $w_i=\bigvee(w_i^1,\dots ,w_i^{q_i})$, for $i=1, 2$, and that $u=\bigvee (u^1,\dots, u^s)$ and $v=\bigvee (v^1,\dots ,v^k)$. 

From Definition \ref{partordertrees}, it is immediate to verify that $v\leq_T u$ implies that $\vert v^1\vert \leq \vert u^1\vert$. 

On the other hand, $t_1/w_1\leq_T u\leq_T t_1\perp w_1$ implies that $l= \vert t_1\vert + \vert w_1^1\vert$. But the condition $t_2/w_2\leq_T v\leq_T t_2\backslash w_2$ implies that $\vert v^1\vert = \vert t_1^1<\vert t_1\vert \leq \vert u^1\vert$. So, we get that $v\not\leq_T u$, and the proof is over.\end{proo}
\medskip

As in the first example, the restriction of the facial order on planar rooted trees to binary planar rooted trees, gives the Tamari order $\leq_{Ta}$ on the set ${\mathcal Y}_n$, for $n\geq 1$ (see \cite{HuaTam}, and \cite{Rea}). Using the results of \cite{LodRon1}, it is easy to see that the family ${\bold{\mathcal Y}}:=\{{\mathcal Y}_n\}_{n\geq 1}$, equipped with $\leq_{Ta}$, is also a dendriform poset.
\bigskip

\subsection{The $\Dy^m$ algebra of $m$-simplexes} 
\medskip

Let ${\bold P}=\{ P_n\}_{n\geq 1}$ be a dendriform poset.

\begin{notation} \label{notnordershuffles} For a positive integer $n$, the ${\mbox {Simp}({\bold P})}_n^m$ be the set of all $m$-simplexes of the partially ordered set $P_n$. That is, 
\begin{equation*}{\mbox {Simp}({\bold P})_n^m}:=\{ (x_1,\dots,x_m)\in P_n^{m}\mid x_1\leq x_2\leq \dots \leq x_m\}\subseteq P_n^{\times m}.\end{equation*}

For $m$ fixed, let $\Ordm({\bold P}):= \bigoplus_{n\geq 1}\K[ {\mbox {Simp}({\bold P})}_n^m]$ be the graded vector space spanned by the set $\displaystyle {\bigcup_{n\geq1} {\mbox {Simp}({\bold P})}_n^m}$ of all $m$-simplexes, for $n\geq 1$.\end{notation} 

We know that the space $V_{\small{\bold P}}$ has a natural structure of dendriform algebra. The aim of the present subsection is to describe a structure of $\Dy^m$-algebra on $\Ordm({\bold P})$.

\begin{definition} \label{mDyckonposet} Let ${\bold P}$ be a dendriform poset. For a pair of elements ${\underline x}=(x_1,\dots, x_m)\in {\mbox {Simp}({\bold P})}_n^m$ and ${\underline y}=(y_1,\dots ,y_m)\in {\mbox {Simp}({\bold P})}_r^m$,, and any $0\leq i\leq m$, let ${\mathcal I}^i({\underline x},{\underline y})$ be the set of all elements ${\underline u}\in {\mbox {Simp}({\bold P})}_{n+r}^m$ satisfying that \begin{enumerate}[(i)]
\item $x_j/y_j\leq u_j\leq x_j\perp y_j$, for $1\leq j\leq {\mbox{\it m-i}}$
\item $x_j\top y_j\leq u_j\leq x_j\backslash y_j$, for ${\mbox{\it m-i}} < j \leq m$.\end{enumerate}\end{definition}
\medskip

In order to define products $*_i$ on the vector space $\Ordm({\bold P})$, for $0\leq i\leq m$, we need to prove some results.

\begin{lemma} \label{lemexmp1} Let $0\leq i<j\leq m$. For any pair of elements $({\underline x},{\underline y})\in {\mbox {Simp}({\bold P})}^m$, we have that the set ${\mathcal I}^i({\underline x},{\underline y})$ is equal to the set of all elements 
${\underline u} \in {\mbox {Simp}({\bold P})}_{n+r}^m$ satisfying that \begin{enumerate}[(i)]
\item $x_k/y_k\leq u_k\leq x_k\backslash y_k$, for $1\leq k\leq m-j$,
\item $x_k/y_k\leq u_k\leq x_k\perp y_k$, for $m-j < k\leq m-i$,
\item $x_j\top y_j\leq u_j\leq x_j\backslash y_j$, for $m-i < k\leq m$.
\end{enumerate}\end{lemma}
\medskip

\begin{proo} As ${\bold P}$ is a dendriform poset, the interval $[x/y; x\backslash y]$ is equal to the disjoint union of $[x/y; x\perp y]$ and $[x\top y; x\backslash y]$, for any elements $x,y\in \bigcup_{n\geq 1}P_n$. 

So, any ${\underline u}\in {\mathcal I}^i({\underline x},{\underline y})$ satisfies that $x_k/y_k\leq u_k\leq x_k\backslash y_k$, for $1\leq k\leq m-j$, which proves the inclusion of ${\mathcal I}^i({\underline x},{\underline y})$.
\medskip

We need to see that any ${\underline u}\in {\mbox {Simp}({\bold P})}_{n+r}^m$, such that \begin{enumerate} [(i)]
\item $x_k/y_k\leq u_k\leq x_k\backslash y_k$, for $1\leq k\leq m-j$
\item $x_k/y_k\leq u_k\leq x_k\perp y_k$, for $m-j <k\leq m-i$,\end{enumerate}
satisfies that $x_k/y_k\leq u_k\leq x_k\perp y_k$, for $1\leq k\leq m-i$. 

As $i <j$, we get that $m-j < m-i$. Suppose that for some $1\leq k\leq m-j$, we have that $x_k\top y_k\leq u_k\leq x_k\backslash y_k$. 

\noindent In this case $u_k\leq u_{m-i}$ and $x_{m-i}/y_{m-i}\leq u_{m-i}\leq x_{m-i}\perp y_{m-i}$, which contradicts condition $(5)$ of Definition \ref{defdendposet}. 

\noindent So, $x_k/y_k\leq u_k\leq x_k\perp y_k$, for all $1\leq k\leq m-j$, which implies that ${\underline u}\in {\mathcal I}^i({\underline x},{\underline y})$.\end{proo}
\medskip

The following Lemma is a consequence of similar arguments than the ones applied in the proof of Lemma \ref{lemexmp1}.

\begin{lemma} \label{lemexmp2} Let $0\leq j\leq i-1 < m$. For any pair of elements $({\underline x},{\underline y})\in {\mbox {Simp}({\bold P})}^m$, we have that the set ${\mathcal I}^i({\underline x},{\underline y})$ is equal to the set of all elements 
${\underline u}\in {\mbox {Simp}({\bold P})}_{n+r}^m$ satisfying that\begin{enumerate}[(i)]
\item $x_k/y_k\leq u_k\leq x_k\backslash y_k$, for $1\leq k\leq m-i$,
\item $x_k/y_k\leq u_k\leq x_k\backslash y_k$, for $m-i <k\leq m-j$,
\item $x_j\top y_j\leq u_j\leq x_j\backslash y_j$, for $m-j < k \leq m$.\end{enumerate}\end{lemma}
\medskip

\begin{definition} \label{mDyckstruct} The product $*_i: \Ordm({\bold P})\otimes \Ordm({\bold P})\longrightarrow \Ordm({\bold P})$, for $0\leq i\leq m$, is defined by
\begin{equation} {\underline x}*_i{\underline y}:=\sum_{{\underline u}\in {\mathcal I}^i({\underline x},{\underline y})}u,\end{equation}
for any pair of elements ${\underline x}$ and ${\underline y}$ in ${\mbox {Simp}({\bold P})}^m$. The products are extended by linearity to the vector space $\Ordm({\bold P})$.\end{definition}
\medskip

\begin{theorem} \label{thdycmexamp} For any dendriform poset ${\bold P}$ and any $m\geq 1$, the vector space $\Ordm({\bold P})$, equipped with the products $*_i$, for $0\leq i\leq m$, is a $\Dy^m$ algebra.\end{theorem}
\medskip

\begin{proo} We need to prove that $(\Ordm({\bold P}),*_0,\dots ,*_m)$ satisfies the relations of Definition \ref{defDyckmalgebra}. For $n=1$, the result was proved in Proposition \ref{dendestposet}.

For $n\geq 2$, we have that
From Definition \ref{defdendposet}, we have that\begin{align}
\sum_{y_k/z_k\leq u_k\leq y_k\perp z_k}(\sum_{x_k/u_k\leq v_k\leq x_k\perp u_k}v_k)& =  \sum_{x_k/y_k\leq w_k\leq x_k\backslash y_k}(\sum_{w_k/z_k\leq v_k\leq w_k\perp z_k}v_k),\\
\sum_{y_k\top z_k\leq u_k\leq y_k\backslash z_k}(\sum_{x_k/u_k\leq v_k\leq x_k\perp u_k}v_k) &= \sum_{x_k/y_k\leq w_k\leq x_k\perp y_k}(\sum_{w_k\top z_k\leq s_k\leq w_k\backslash z_k}s_k),\\
\sum_{y_k/ z_k\leq u_k\leq y_k\backslash z_k}(\sum_{x_k\top u_k\leq v_k\leq x_k\backslash u_k}v_k)& = \sum_{x_k\top y_k\leq w_k\leq x_k\backslash y_k}(\sum_{w_k\top z_k\leq v_k\leq w_k\backslash z_k}v_k),\end{align}
for any $1\leq k\leq m$ and ${\underline x}, {\underline y}, {\underline z}\in {\mbox {Simp}({\bold P})}^m$.

In all the equalities above, from condition $(4)$ of Definition \ref{defdendposet}, we have that if $v_1\leq \dots \leq v_m$, then $u_1\leq \dots \leq u_m$ and $w_1\leq \dots \leq w_m$. 
\begin{enumerate} \item Let $0\leq i<j\leq m$. We have that
\begin{equation*}{\underline x}*_i({\underline y}*_j{\underline z})=\sum_{v_1\leq \dots \leq v_m} {\underline v},\end{equation*}
where the sum is taken over all the elements ${\underline v}\in {\mbox {Simp}({\bold P})}^m$ such that there exists ${\underline u}\in {\mbox {Simp}({\bold P})}^m$ satisfying that \begin{enumerate}
\item $y_k/z_k\leq u_k\leq y_k\perp z_k$ and $x_k/u_k\leq v_k\leq x_k\perp u_k$, for $1\leq k\leq m-j$,
\item $y_k\top z_k\leq u_k\leq y_k\backslash z_k$ and $x_k/u_k\leq v_k\leq x_k\perp u_k$, for $m-j < k\leq m-i$, 
\item $y_k\top z_k\leq u_k\leq y_k\backslash z_k$ and $x_k\top u_k\leq v_k\leq x_k\backslash u_k$, for $m-i < k\leq m$.\end{enumerate}

As $m-j < m-i$ and $u_1\leq \dots \leq u_m$, Lemma \ref{lemexmp2} implies that we may change condition$(c)$ in the last paragraph, for the following condition

\noindent $(c\rq )$ for $m-i < k\leq m$, $y_k/ z_k\leq u_k\leq y_k\backslash z_k$ and $x_k\top u_k\leq v_k\leq x_k\backslash u_k$.
\medskip

Applying conditions $(1.1)$, $(1.2)$ and $(1.3)$ above, we get that 
\begin{equation*}{\underline x}*_i({\underline y}*_j{\underline z}) =\sum_{v_1\leq \dots \leq v_m} {\underline v},\end{equation*}
where the sum is taken over all the elements ${\underline v}\in {\mbox {Simp}({\bold P})}^m$ such that there exists ${\underline w}\in {\mbox {Simp}({\bold P})}^m$ satisfying that\begin{enumerate}
\item $x_k/y_k\leq w_k\leq x_k\backslash y_k$ and $w_k/z_k\leq v_k\leq w_k\perp z_k$, for $1\leq k\leq m-j$, 
\item  $x_k/y_k\leq w_k\leq x_k\perp y_k$ and $w_k\top z_k\leq s_k\leq w_k\backslash z_k$, for $m-j < k\leq m-i$,
\item  $x_k\top y_k\leq w_k\leq x_k\backslash y_k$ and $w_k\top z_k\leq v_k\leq w_k\backslash z_k$, for $m-i < k\leq m$.\end{enumerate}
\medskip

Again, as $m-j < m-i$ and $w_1\leq \dots \leq w_m$,  Lemma \ref{lemexmp1} implies that condition $(a)$ may be changed into

\noindent $(a\rq)$ $x_k/y_k\leq w_k\leq x_k\perp y_k$ and $w_k/z_k\leq v_k\leq w_k\perp z_k$, for $1\leq k\leq m-j$,

\noindent which implies that ${\underline x}*_i({\underline y}*_j{\underline z}) = ({\underline x}*_i{\underline y})*_j{\underline z}$.
\medskip

\item Let $0\leq i\leq m$, it is easy to show that\begin{enumerate}
\item $\sum_{j=0}^i {\underline y}*_j{\underline z} = \sum_{u_1\leq\dots \leq u_m}{\underline  u}$, 
where the sum is taken over all the elements ${\underline u}$ satisfying that\begin{enumerate}
\item $y_k/z_k\leq u_k\leq y_k\perp z_k$, for $1\leq k\leq m-i$,
\item $y_k/ z_k\leq u_k\leq y_k\backslash z_k$, for $m-i < k\leq m$,\end{enumerate}
\item $\sum_{j=i}^m {\underline x}*_j{\underline y} = \sum_{w_1\leq \dots \leq w_m}{\underline w}$, 
where the sum is taken over all the elements ${\underline w}$ satisfying that \begin{enumerate}
\item $x_k/y_k\leq w_k\leq x_k\backslash y_k$, for $1\leq k\leq m-i$,
\item $x_k\top y_k\leq w_k\leq x_k\backslash y_k$, for $m-i < k\leq m$,\end{enumerate} 
for elements ${\underline x},  {\underline y}$ and  ${\underline z}$ in ${\mbox {Simp}({\bold P})}^m$.
\end{enumerate}
\medskip

We get that \begin{equation*} {\underline x}*_i \sum_{j=0}^i {\underline y}*_j{\underline z}=\sum_{v_1\leq \dots \leq v_m} {\underline v},\end{equation*} where the sum is taken over all elements ${\underline v}\in {\mbox {Simp}({\bold P})}^m$ such that there exists ${\underline u}$ satisfying that \begin{enumerate}
\item $y_k/z_k\leq u_k\leq y_k\perp z_k$ and $x_k/u_k\leq v_k\leq x_k\perp u_k$, for $1\leq k\leq m-i$,
\item  $y_k/ z_k\leq u_k\leq y_k\backslash z_k$ and $x_k\top u_k\leq v_k\leq x_k\backslash u_k$, for $m-i < k\leq m$.\end{enumerate}

Note that $v_1\leq \dots \leq v_m$, implies that ${\underline u}\in {\mbox {Simp}({\bold P})}^m$ from condition $(4)$ of Definition \ref{defdendposet}.
\medskip

Applying conditions $(1.1)$, $(1.2)$ and $(1.3)$, we obtain \begin{equation*} {\underline x}*_i \sum_{j=0}^i {\underline y}*_j{\underline z} =\sum_{v_1\leq \dots \leq v_m} {\underline v},\end{equation*} where\begin{enumerate}[(i)]
\item $x_k/y_k\leq w_k\leq x_k\backslash y_k$ and $w_k/z_k\leq v_k\leq w_k\perp z_k$, for $1\leq k\leq m-i$,
\item $x_k\top y_k\leq w_k\leq x_k\backslash y_k$ and $w_k\top z_k\leq v_k\leq w_k\backslash z_k$, for $m-i < k\leq m$,\end{enumerate}
for some ${\underline w}\in {\mbox {Simp}({\bold P})}^m$.
\medskip

As \begin{equation*}\sum_{j=i}^m{\underline x}*_j{\underline y} =\sum_{w_1\leq \dots \leq w_m}{\underline w},\end{equation*} 
with \begin{itemize}
\item $x_k/y_k\leq w_k\leq x_k\backslash y_k$, for $1\leq k\leq m-i$, 
\item $x_k\top y_k\leq w_k\leq x_k\backslash y_k$, for $m-i < k\leq m$. \end{itemize}
\medskip

We get that \begin{equation*}{\underline x}*_i(\sum_{j=0}^i {\underline y}*_j{\underline z})= (\sum_{j=i}^m {\underline x}*_j{\underline y})*_i {\underline z},\end{equation*} which ends the proof.
\end{enumerate}\end{proo}

\bigskip

\section{$\Dy^m$ algebras and $m${-}Dyck paths} 
\medskip

The dimension of the subspace of homogeneous elements of degree $n$ of the free $\Dy^m$ algebra ${\mathfrak D}^m$ is the number of $m$-Dyck paths of size $n$. We define a $\Dy^m$ algebra structure (isomorphic to ${\mathfrak D}^m$) on the vector space spanned by the set of all $m$-Dyck paths, and prove that this structure may be described in terms of the $m$-Tamari lattice introduced by F. Bergeron in \cite{Ber}.

\subsection{$m${-}Dyck paths}

We describe basic notions bout $m$-Dyck paths. For more detailed constructions and the proofs of the results we refer to \cite{BerPre}, \cite{BouFusPre} and \cite{BMChPR2013}.

\begin{definition} \label{Dyckpath} For $m, n\geq 1$, an {\it $m$-Dyck path of size $n$} is a path on the real plan 
${\mathbb R}^2$, starting at $(0,0)$ and ending at $(2nm,0)$, consisting of up steps $(m,m)$ and down steps $(1,-1)$, 
which never goes below the $x$-axis. Note that the initial and terminal points of each step lean on $\ZZ_+^2$.\end{definition}

\begin{notation} \label{rho} We denote by $\Dyc _n^m$ the set of all $m$-Dyck paths of size $n$. Define $\Dyc _0^m :=\{ \bullet\}$, for $m\geq 1$.

 For any $m\geq 0$, we denote by $\rho_m\in \Dyc _1^m$ the unique $m$-Dyck path of size one.\end{notation}

In order to define constructions on Dyck paths, we use the notation employed by M. Bousquet-M\'elou, E. Fusy and the second author in \cite{BouFusPre}.

\begin{notation} \label{updown} Let $P$ be an $m$-Dyck path. We denote by $\up(P)$ the set of up steps of $P$, and by $\dw(P)$ the set of down steps of $P$.\end{notation}

\begin{definition} \label{excursion} Let $u\in \up(P)$ be an up step of an $m$-Dyck path $P$, 
the {\it rank} of $u$ is $k$ if $u$ is the $k^{\rm th}$ up step of $P$, counting from left to right. 

A down step $d$ is at {\it level} $k$ if the last up step $u$ preceding $d$ has rank $k$. \end{definition}

As any up step $u$ in a Dyck path $P$ is determined by its rank, from now on we identify them, 
and denote the set of up steps of a Dyck path of size $n$ as $\up(P) =\{ 1, \dots ,n\}$.

\begin{example} The down steps of the following $2$-Dyck path are colored with their levels

\begin{figure}[h]
\includegraphics[scale=0.6]{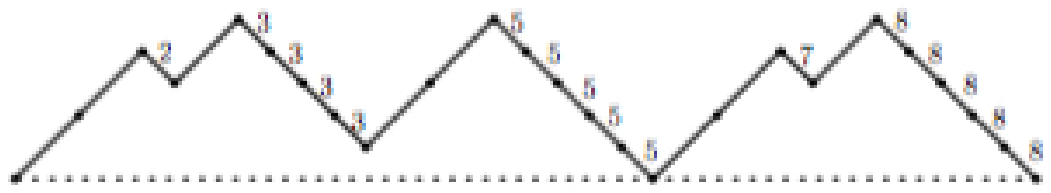}
\end{figure}
\end{example}

\begin{notation}  \label{not:Ln(P)} For a path $P\in \Dyc_n^m$ and an integer $1\leq k\leq n$, we denote by $\dw _k(P)$ the set of down steps of level $k$ of $P$ and by  $L_k(P)$ the number of elements of $\dw _k(P)$. Note that $\dw _k(P)$ may be the empty set, and in this case $L_k(P)=0$. When no confusion is possible, we shall denote the last term of the sequence $L_n(P)$ simply by $L(P)$.

Note that $0\leq \L_1(P)+\dots +L_j(P)\leq mj$, for $1\leq j \leq n$. As a Dyck path $P\in \Dyc_n^m$ is uniquely determined by the sequence $L_1(P),\dots , L_{n-1}(P),L(P)$, we denote $P= ((L_1(P),\dots , L_(P)))$.

For any integer $1\leq k\leq n$ such that $\dw_k(P)\neq \emptyset$, we identify the set $\dw _k(P)$ with the sequence of its down steps $\dw _k(P)= d_1^kd_2^k\dots d_{L_k(P)}^k$, ordered from left to right, and from top to bottom.\end{notation}
\bigskip

\subsection{Basic operations on Dyck paths}
\medskip

We want to describe basic operations on Dyck paths that we need in the sequel.

\begin{definition} \label{suspension-concatenation}
Let $P$ and $Q$ be two $m$-Dyck paths of sizes $n$ and $r$, respectively. For $0\leq i\leq L(P)$, 
define the $i^{th}$-{\it concatenation} of $P$ and $Q$, denoted $P\t_i Q$, as the Dyck path of size $n+r$ given by
\begin{equation*} P\t_i Q:=((L_1(P),\dots ,L_{n-1}(P), L(P)-i,L_1(Q),\dots ,L_{r-1}(Q), L(Q)+i)).\end{equation*}\end{definition}

\begin{definition} \label{defprime} An $m$-Dyck path $P$ is called {\it prime} if there does not exist a pair of $m$-Dyck paths $Q$ and $R$ such that $P = Q\t_0 R$, with $P\neq Q$ and $P\neq R$.\end{definition}

\begin{remark} \label{remprime} For any $m$-Dyck path of size $n$ there exist a unique composition $(n_1,\dots ,n_r)$ of $n$ (with $n_i\geq 1$ for each $i$) and a unique family of prime Dyck paths $P_1\in \Dyc_{n_1}^m,\ \dots , P_r\in \Dyc_{n_r}^m$
 such that $P = P_1\t_0 \dots \t_0 P_r$.
\end{remark}
\medskip

The proof of the following Lemma is immediate.

\begin{lemma} \label{lemma1} Let $P\in \Dyc_{n_1}^m$ be a prime Dyck path and let $Q\in \Dyc_{n_2}^m$ be another Dyck path. For any $1\leq j\leq L(P)$, the Dyck path $P\t _j Q$ is prime.\end{lemma}
\medskip

\begin{definition}\label{defcoloring} Let $P$ be an $m$-Dyck path of size $n$. The {\it standard coloring} of $P$ is a map $\alpha _P$ from the set of down steps $\dw(P)$ to the set $\{1,\dots ,n\}$, described recursively as follows\begin{enumerate}
\item For $P = \rho _m\in \Dyc_1^m$, $\alpha _{\rho _m}$ is the constant function $1$.
\item For $P = \bigvee_d (P_0,\dots ,P_m)$, with $P_j\in \Dyc_{n_j}^m$, the set of down steps of $P$ is the disjoint union 
\begin{equation*}\dw (P) = \{ 1,\dots ,m\}\bigcup \dw(P_0) \bigcup \dots \bigcup  \dw(P_m),\end{equation*}
where the first subset $\{1,\dots ,m\}$ corresponds to the down steps of $\rho _m$. 

The map $\alpha _P$ is defined by
\begin{equation*}\alpha _P(d) = \begin{cases} 1,&{\rm for}\ d\in \{1,\dots ,m\},\\
\alpha _{P_j}(d)+n_0+\dots +n_{j-1}+1,&{\rm for}\ d\in \dw(P_j),\end{cases}\end{equation*}
where $0\leq j\leq m$.
\end{enumerate}\end{definition}
\medskip

For instance, the standard coloring of the $2$-Dyck path $P=((0,1,4, 0, 5, 0,1,5))$ is

\begin{figure}[h]
\includegraphics[scale=0.6]{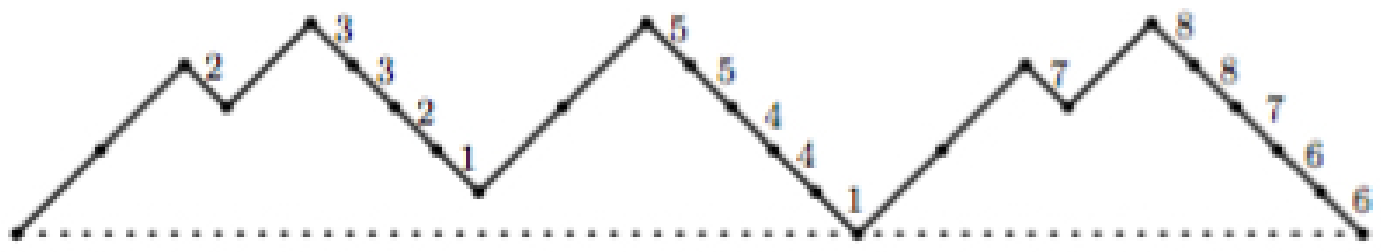}
\end{figure}

\begin{remark} \label{remcoloring} Let $P$ be an $m$-Dyck path of size $n$.
\begin{enumerate} 
\item We have that $\vert \alpha _P^{-1}(i)\vert = m$, for any $1\leq i\leq n$.
\item Let $\dw^k(P) = d_{1}^k,d_2,\dots ,d_{L_k(P)}^k$ be the ordered sequence of down steps of level $k$ of $P$, for $1\leq k\leq n$. 

\noindent The word $\omega _k^P= \alpha _P(d_{1}^k)\dots  \alpha _P(d_{L_k(P)}^k)$ 
is decreasing for the usual order of the natural numbers. Moreover, the first $m$ digits of $\omega _n^P$ are $n$\rq s. 
\item If $Q$ is another $m$-Dyck path, then $\dw (P\t _i Q) = \dw (P)\bigcup \dw (Q)$, and $\alpha _{P\t_iQ}$ is described by
\begin{equation*}\alpha _{P\t_iQ}(d)= \begin{cases} \alpha _P(d),&{\rm for\ any}\ d\ {\rm which\ belongs\ initially\ to}\ P,\\
\alpha _Q(d) + n,&{\rm for\ any}\ d\ {\rm which\ belongs\ initially\ to}\ Q,\end{cases}\end{equation*}
for any $0\leq i\leq L(P)$.
\end{enumerate}
\end{remark}
\medskip

\subsection{The $\Dy^m$ algebra structure on $m$-Dyck paths}
\medskip

\begin{definition} \label{defpartition} For any positive integer $n$, a {\it weak composition} of $n$ with $r+1$ parts is an ordered collection 
of non-negative integers ${\underline {\lambda}}= (\lambda_0,\dots ,\lambda _r)$ such that $\lambda _0+\dots +\lambda_r = n$. We say that the {\it length} of $\lam$ is $r+1$.\end{definition}

\begin{notation} \label{notLambda} The set of all weak compositions of $n$ with $r+1$ parts is denoted by $\Lambda _r^n$. Given an $m$-Dyck path $P$ of size $n$ we denote the set of all weak compositions of $L(P)$ of length $r+1$ by 
$\Lambda _r(P)$.\end{notation}

Let $Q\in \Dyc _{s}^m$ be an $m$-Dyck path, such that $Q = Q_1\t_0 \dots \t_0 Q_r$ with $Q_j\in \Dy^m$ prime, for $1\leq j\leq r$. Given another Dyck path $P\in \Dyc_{n}^m$, suppose that $\lam =(\lambda _0,\dots ,\lambda _r)$ is a weak 
composition of $L(P)$. 
Define a Dyck path $P*_{\lam} Q$ of size $n+s$ by the formula
\begin{equation} P*_{\lam} Q := ((((P\t _{\lambda _1+\dots +\lambda _r}Q_1)\t _{\lambda _2+\dots + \lambda _r} Q_2)\t _{\lambda _3+\dots +\lambda _r}\dots )\t _{\lambda _r}Q_r).\end{equation}

The product $*_{\lam}$ just divides the ordered set $\dw_n(P)$ of down steps of level $n$ of $P$ and glue, in order, the $i^{th}$ piece at the end of the path $Q_i$. If $\lambda _0 > 0$, the first $\lambda _0$ steps of $\dw_n(P)$ remain at the end 
of $P$.

\begin{example} Let $P = ((2, 3, 1, 6))\in \Dyc _4^3$ and let $Q = ((1,4,4,3, 2,3,4))$ be a $3$-Dyck path of size $7$, where we denote the Dyck paths following Notation 
\ref{not:Ln(P)}. We have that that $Q = ((1,4,4))\t _0 \rho_3\t_0 ((2,3,4))$. 

The word on the top level of $P$ is $\omega _4^P:= 444331$. Consider the weak composition $\lam = (1,2,2,1)$ of $L(P)= 6$ of length $4$. 
The path $P\t_{(1,2,2,1)}Q$ is the following

\begin{figure}[h]
\includegraphics[scale=0.6]{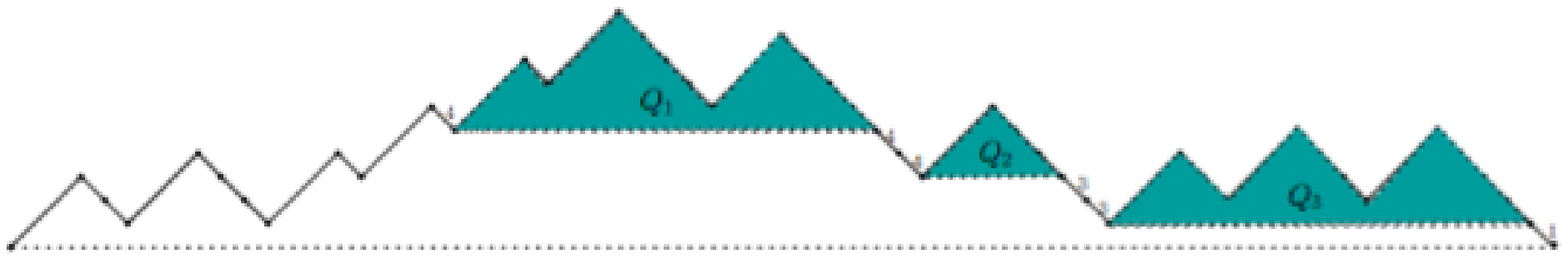}
\end{figure}
\end{example}

The last point of Remark \ref{remcoloring} implies that for any $P\in \Dyc_n^m$, any $Q = Q_1\t_0\dots \t_0 Q_r$, with $Q_i$ prime for $1\leq i\leq r$, and any $\lam \in \Lambda _r(P)$, the set of down steps of 
$P*_{\lam}Q$ is
\begin{equation*} \dw (P*_{\lam}Q) = \dw (P)\ \bigcup\ \dw (Q),\end{equation*}
and the standard coloring $\alpha _{P*_{\lam}Q}$ is described by
\begin{equation*}\alpha _{P*_{\lam}Q}(d) = \begin{cases} \alpha _P(d),& {\rm for}\ d\in \dw(P),\\
\alpha _Q(d) + n,& {\rm for}\ d\in \dw(Q).\end{cases} \end{equation*}
\medskip

\begin{notation} \label{notparti} Let $P$ be a Dyck path of size $n$, with $\dw _{n}(P)=d_{1}^n,\dots ,d_{L(P)}^n$, and let $\lam = (\lambda _0,\dots ,\lambda _r)$ be a weak composition of $L(P)$. For $0\leq i\leq m$, we denote by 
$\Lambda _r ^i(P)$ the set of all weak compositions $\lam$ of length $r+1$ such that the restriction $\alpha_P (d_{L(P){-}\lambda _r+1}^n),\dots ,\alpha_P(d_{L(P)}^n)$ of the word $\omega _{n}^P$ to its last $\lambda_r$ letters satisfies the following conditions:\begin{enumerate}
\item any digit in the word $\alpha_P (d_{L(P)-\lambda _r+1}^n),\dots ,\alpha_P(d_{L(P)}^n)$ appears at most $i$ times,
\item there exists at least one integer $1\leq i_0\leq n$ such that $i_0$ appears exactly $i$ times in $\alpha_P (d_{L(P){-}\lambda _r+1}^n),\dots ,\alpha_P(d_{L(P)}^n)$.\end{enumerate}\end{notation}
\medskip

For example, for $P = ((0,2,1,3,4))\in \Dy c_5^2$, we get that $\omega _{5}^P= 5,5,1,1$.
So, $\lam _1= (1, 1, 2)$ belongs to $\Lambda _2^2(P)$, while $\lam _2 = (0, 3,1)$ belongs to $\Lambda _2^1(P)$.

Observe that
\begin{equation*}\Lambda _r^0(P) = \{ (\lambda _0,\dots ,\lambda _{r-1}, 0)\mid \lambda _0+\dots +\lambda_{r-1} = L(P)\ {\rm and}\ r\geq 1\}.\end{equation*}

The set of all weak compositions of $L(P)$ is the disjoint union ${\displaystyle \bigcup _{r\geq 0}\bigl (\bigcup _{i=0}^m\Lambda _r^i (P)\bigr )}$, for any $m$-Dyck path $P$ of size $n$. 
 \medskip
 
 The following result is a straightforward consequence of Lemma \ref{lemma1} and the definition of $*_{\lambda}$.

\begin{lemma} \label{lemma2} Let $P = P_1\t _0\dots \t_0 P_s$ in $\Dyc _n^m$ and $Q= Q_1\t_0\dots \t_0Q_r$ in $\Dyc _q^m$ be two Dyck paths, where $P_1,\dots ,P_s,Q_1,\dots ,Q_r$ are prime Dyck paths, and let 
$\lam \in \Lambda _r^i(P)$ be a weak composition. We have that\begin{enumerate}
\item if $i > 0$, then
\begin{equation*}P*_{\lam}Q = P_1\t _0\dots \t_0 P_{s-1}\t_0 (P_s*_{\lam } Q),\end{equation*}
where $P_s*_{\lam } Q$ is prime.
\item if $i=0$, then $\lam = (\lambda _0,\dots ,\lambda _{r-1},0)$ and
\begin{equation*} P*_{\lam}Q = P_1\t _0\dots \t_0 P_{s-1}\t_0 (P_s*_{\lam } (Q_1\t _0 \dots \t_0 Q_{j_0}))\t _0 Q_{j_0+1}\t_0\dots \t_0 Q_r,\end{equation*}
where $j_0$ is the maximal element of $\{0,\dots , r{-}1\}$ such that $\lambda _{j_0}\neq 0$.
\end{enumerate}
\end{lemma}
\medskip
 
The product on the graded vector space $\K[\Dyc^m]$, spanned by the set of all $m$-Dyck paths, is defined as follows.
 
 \begin{definition} \label{defprodi} Let $P\in \Dyc _n^m$ and $Q\in \Dyc _s^m$ be two Dyck paths, such that $Q = Q_1\t _0\dots \t_0 Q_r$ with $Q_i$ prime, $1\leq i\leq r$. For any integer $0\leq j\leq m$, define
\begin{equation} P*_jQ= \sum _{\lam\in \Lambda _r^j(P)} P*_{\lam}Q.\end{equation}
The product extends in a unique way to a linear map from $\K[\Dyc^m]\ot \K[\Dyc^m]$ to $\K[\Dyc^m]$.\end{definition}
\medskip

\begin{example} Consider the $2$-Dyck paths $P = ((1,3))\in \Dyc_2^2$ and $Q = ((0, 2, 4,2)) = ((0,2,4))\t_0 \rho_2\in \Dyc_4^2$.
Computing the products $P*_0Q$ and $P*_1Q$ , we get that
\begin{align*}  P*_0Q &= P*_{(3,0,0)}Q + P*_{(2,1,0)}Q + P*_{(1,2,0)}Q + P*_{(0,3,0)}Q =\\
&((1,3,0,2,4,2)) + ((1,2,0,2,5,2)) + ((1,1,0,2,6,2)) + ((1,0,0,2,7,2)) =\end{align*}
\begin{figure}[h]
\includegraphics[scale=0.6]{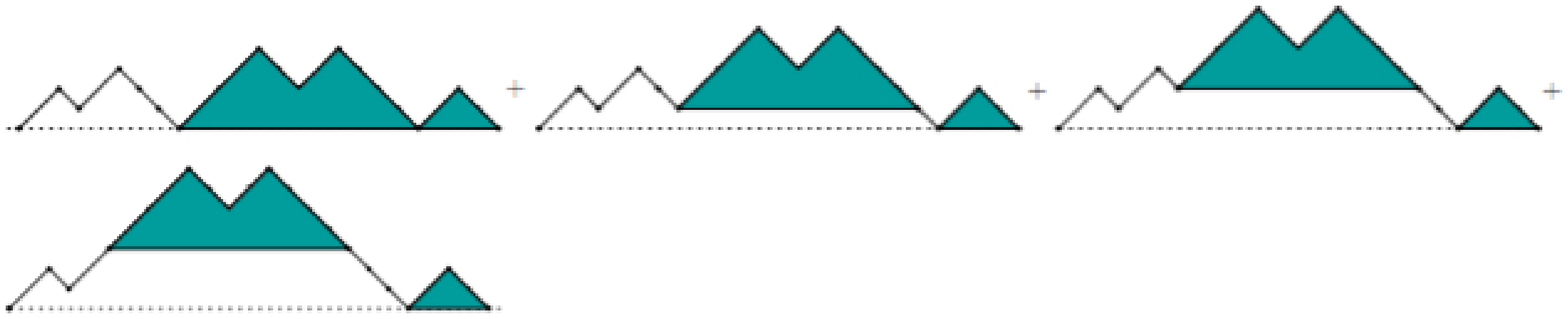}
\end{figure}
\end{example}

\begin{proposition} \label{proprelac} Let $P\in \Dyc _n^m$ and $Q = Q_1\t _0\dots \t _0 Q_r\in \Dyc _p^m$ be two Dyck paths, with $Q_j\in \Dyc _{p_j}^m$ prime for $1\leq j\leq r$. \begin{enumerate}
\item For nonnegative integers $s\geq 1$ and $0\leq i < j\leq m$, the map 
\begin{equation*}\psi_{ij}(P,Q):\{ \Lambda _r^i(P)\times \Lambda _s^j(Q) \longrightarrow \{ (\lam ,{\underline {\delta}})\mid \lam \in \Lambda_r^i(P)\ {\rm and}\ {\underline {\delta}}\in \Lambda _s^j(P*_{\lam}Q)\},\end{equation*}
which sends $(\lam , {\underline {\tau}})\mapsto (\lam , {\underline{\delta}} :=(\tau_0,\dots ,\tau _{s-1},\tau_s +\lambda_r))$ is bijective.
\item For any integer $0\leq i\leq m$, the map $\psi _i^1(P,Q)(\lam ,{\underline {\tau}}): =$
\begin{equation*}((\lambda _0,\dots ,\lambda _{r-1},\lambda _r + \dots + \lambda _{r+s{-}j_{\underline {\tau}}}) , (\tau_0,\dots ,\tau_{j_{\underline {\tau}}}+\lambda _r,\lambda _{r+1},\dots ,\lambda _{r+s{-}j_{\underline {\tau}}})),\end{equation*} defines a bijection from the set $\{(\lam ,{\underline {\tau}})\mid {\underline {\tau}}\in \Lambda _s^0(Q)\ {\rm and}\ \lam \in \Lambda _{r+s{-}j_{\underline {\tau}}}^i (P)\}$ to the set 
\begin{equation*}\{ ({\underline {\gamma }}, {\underline {\delta }})\mid {\underline {\gamma }} \in {\displaystyle \bigcup_{j=i}^m\Lambda _r^j(P)}\ {\rm and}\  {\underline {\delta }}\in \Lambda _s^i(P*_{\underline {\gamma }}Q)\ {\rm such\ that}\ \delta _s \leq \gamma _r\},\end{equation*}  
where $j_{\underline {\tau}}$ is the maximal integer $0\leq j\leq s-1$ such that $\tau _j > 0$, and $\bigcup $ denotes the disjoint union. 
\item For any integer $0\leq i\leq m$, the map
\begin{equation*}\psi _i^2(P,Q):\Lambda _r^i(P)\t ({\displaystyle \bigcup _{j=1}^i \Lambda _s^j(Q)})\longrightarrow \{ ({\underline{\gamma}} , {\underline {\delta}})\mid {\underline{\gamma}} \in \Lambda _r^i(P)\ {\rm and}\ {\underline {\delta}}\in \Lambda _s^i(P*_{\underline {\gamma}}Q)\ {\rm such\ that}\ \gamma _r < \delta _s\},\end{equation*}
which maps $(\lam , {\underline {\tau}})\mapsto (\lam , {\underline{\delta}}:= (\tau_0,\dots ,\tau _{s-1}, \tau _s+\lambda _r)$, is bijective.
\end{enumerate}
\end{proposition}
\medskip

\begin{proo} $(1)$ For the first point, let $\lam \in \Lambda _r^i(P)$ and ${\underline {\tau}}\in  \Lambda _s^j(Q)$ be two weak compositions. 

If $\dw_n(P) = d_{1}^n,\dots ,d_{L(P)}^n$ and $ \dw_p(Q)= d_{1}^p,\dots ,d_{L(Q)}^p$, then
\begin{equation*}\dw _{n+p}(P*_{\lam}Q) = d_{1}^p,\dots ,d_{L(Q)}^p, d_{L(P){-}\lambda_r+1}^n,\dots ,d_{L(P)}^n.\end{equation*} 

The map $\psi _{ij}$ is defined by the formula
\begin{equation*}\psi _{ij}(\lam, {\underline {\tau}}) := (\lam , {\underline{\delta}} :=(\tau_0,\dots ,\tau _{s{-}1},\tau_s +\lambda_r)).\end{equation*} 

Clearly, $\lam $ belongs to $ \Lambda_r^i(P)$. On the other hand, 
\begin{equation*}\dw _{n+p}(P*_{\lam}Q) = d_1^p,\dots ,d_{L(Q)}^p,d_{L(P){-}\lambda _r+1}^n,\dots ,d_{L(P)}^n,\end{equation*} 
 which implies that the subset of the last 
$\tau_s +\lambda_r$ down steps of $P*_{\lam}Q$ is
$d_{L(Q){-}\tau _s+1}^p, \dots ,d_{L(Q)}^p, d_{L(P){-}\lambda _r+1}^n,\dots ,d_{L(P)}^n$.
\medskip

Note that \begin{enumerate}[(a)]
\item $\alpha _{P*_{\lam}Q}(d_{L(P){-}\lambda _r+1}^n)\dots \alpha_{P*_{\lam}Q}(d_{L(P)}^n)$ is a sequence of elements in the set $\{1,\dots n\}$ such that any digit appears at most $i$ times.
\item $\alpha _{P*_{\lam}Q}(d_{L(Q){-}\tau _s+1}^p)\dots \alpha _{P*_{\lam}Q}(d_{L(Q)}^p)$ is a sequence of elements in the set $\{n+1,\dots ,n+p\}$ where there exists at least one digit that appears $j$ times, and no digit appears more than $j$ times.\end{enumerate}
So, ${\underline{\delta}}$ belongs to $\Lambda _s^j(P*_{\lam}Q)$.
\medskip

For any pair of weak compositions $\lam \in \Lambda _r^i(P)$ and ${\underline {\delta}}\in \Lambda _s^j(P*_{\lam}Q)$, we get that $\omega_{n+p}^{P*_{\lam}Q}$ is equal to
\begin{equation*}\alpha _Q(d_{L(Q){-}\delta_s+1}^p)+n,\dots ,\alpha _Q(d_{L(Q)}^p)+n,\alpha _P(d_{L(P){-}\lambda _r+1}^n),\dots ,\alpha _P(d_{L(P)}^n).\end{equation*}

As the expression $\alpha_P(d_{L(P){-}\lambda_r+1}^n)\dots \alpha _P(d_{L(P)}^n)$ is a word in the alphabet $\{ 1,\dots ,n\}$ such that no digit appears more than $i$ times, and $i < j$, then ${\underline {\tau}} := (\delta _0,\dots ,\delta _s{-}\lambda _r)$ must belong to $\Lambda _s^j(Q)$. 
 \medskip
 
It is immediate to prove that the map $(\lam ,{\underline {\delta}})\mapsto (\lam ,{\underline {\tau}})$ is inverse to $\psi_{ij}(P,Q)$, which ends the proof of $(1)$.
\medskip

$(2)$ If $\lam\in \Lambda _{r+s{-}j_{\underline {\tau}}}^i (P)$ and ${\underline {\tau}}\in \Lambda _s^0(Q)$, then it is immediate to verify that \begin{enumerate}[(i)]
\item ${\underline {\gamma }} = (\lambda _0,\dots ,\lambda _{r-1},\lambda _r + \dots + \lambda _{r+s{-}j_{\underline {\tau}}})$ belongs to $\Lambda _r^j(P)$, for $i\leq j\leq m$,
\item $ {\underline {\delta }} = (\tau_0,\dots ,\tau_{j_{\underline {\tau}}}+\lambda _r,\lambda _{r+1},\dots ,\lambda _{r+s{-}j_{\underline {\tau}}})$ belongs to $\Lambda _s^i(P*_{\underline {\gamma }}Q)$,
\item $\delta _s = \lambda _{r+s{-}j_{\underline {\tau}}}\leq \gamma _r = \lambda _r + \dots + \lambda _{r+s{-}j_{\underline {\tau}}}$.\end{enumerate}

Assume that we have two weak compositions ${\underline {\gamma }} = (\gamma _0,\dots ,\gamma _r)\in {\displaystyle \bigcup _{j=i}^m\Lambda _r^j(P)}$ and $ {\underline {\delta }} = (\delta _0, \dots ,\delta _s) \in \Lambda _s^i(P*_{\underline {\gamma }}Q)$, such that $\delta _s \leq \gamma _r$.
\medskip

Let $j_0$ be the maximal integer $0\leq j_0\leq s{-}1$, such that $\delta _{j_0}+\dots +\delta _s > \gamma _r$. 

Define
\begin{enumerate}[(a)]
\item $\lam := (\gamma _0,\dots , \gamma_{r-1}, \gamma_r{-}\delta_{j_0+1}{-}\dots {-}\delta_s , \delta_{j_0+1},\dots, \delta_s)$,
\item ${\underline {\tau}}:= (\delta _0,\dots ,\delta_{j_0{-}1}, \delta _{j_0}+\dots +\delta _s -\gamma _r, 0,\dots ,0)$.\end{enumerate}
It is clear that $\lam \in \Lambda _{r+s-j_0}^i(P)$, 
${\underline {\tau}}\in \Lambda _s^0(Q)$ and $\psi _i^1(P,Q)(\lam ,{\underline {\tau}}) = ({\underline {\gamma }},  {\underline {\delta }})$, which shows that $\psi _i^1$ is bijective, ending the proof of $(2)$. 
\bigskip

$(3)$ For $\lam \in \Lambda _r^i(P)$ and ${\underline {\tau}}\in \Lambda _s^j(Q)$, for $1\leq j\leq i$, we have that the weak composition $\psi _i^2(P,Q) (\lam , {\underline {\tau}}) = ({\underline{\gamma}}, {\underline{\delta}})$ 
satisfies the following conditions\begin{enumerate}[(i)]
\item ${\underline{\gamma}} = \lam$ belongs to $\Lambda _r^i(P)$,
\item the weak composition ${\underline{\tau}}$ belongs to $\Lambda _s^j(Q)$, for some $1\leq j\leq i$. So, the sequence $ \alpha _Q(d_{L(Q){-}\tau _s+1}^p)+n\dots \alpha _Q(d_{L(Q)}^p)+n$ is a word in the digits of $\{n+1,\dots ,n+p\}$ such that each sequence appears at most $j$ times.

On the other hand, the sequence $\alpha_P(d_{L(P)-\gamma _r+1}^n)\dots \alpha _P(d_{L(P)}^n)$ 
is a word in $\{1,\dots ,n\}$ such that some digit appears exactly $i$ times in it and no digit appears more than $i$ times. 

The sequence $\omega _{n+p}^{P*_{\underline {\gamma}}Q}$, of level $n+p$ of $P*_{\underline {\gamma}}Q$, is equal to
\begin{equation*}\alpha _Q(d_{L(Q){-}\tau _s+1}^p)+n\dots \alpha _Q(d_{L(Q)}^p)+n\alpha_P(d_{L(P)-\gamma _r+1}^n)
\dots \alpha _P(d_{L(P)}^n),\end{equation*}
which shows that ${\underline {\delta}} = (\tau _0,\dots ,\tau_{s-1},\tau _s+\lambda _r)$ belongs to $\Lambda _s^i(P*_{\gamma}Q)$.

\item As ${\underline {\gamma }} = \lam$ and ${\underline {\delta}} = (\tau_0,\dots ,\tau_{s-1}, \tau_s + \lambda _r)$, with $\tau _s > 0$, we get that $\gamma _r < \delta _s$.\end{enumerate}

The map $({\underline {\gamma}}, {\underline {\delta}})\mapsto ({\underline {\gamma}}, (\delta_0,\dots ,\delta _{s-1}, \delta _s- \gamma _r))$ is the inverse map of $\psi _i^2(P,Q)$.
\end{proo}
\medskip

\begin{theorem} \label{theorem1} The binary operations $*_0,\dots , *_m$ defined on $\K[\Dyc^m]$ satisfy the following relations\begin{enumerate}
\item $x*_i (y*_j z) = (x*_iy)*_jz$, for $0\leq i < j\leq m$,
\item $x*_i (y*_0z + \dots + y*_iz) = (x*_iy + \dots +x*_my)*_iz$, for $0\leq i\leq m$,\end{enumerate}
where $x,y,z$ are arbitrary elements of $\K[\Dyc^m]$.\end{theorem}
\medskip

\begin{proo} Clearly, it suffices to prove the relations for any Dyck paths $P$, $Q$ and $Z$. Suppose that $P\in \Dyc _n^m$, $Q = Q_1\t_0\dots \t_0Q_r\in \Dyc _p^m$ and $Z = Z_1\t_0\dots \t_0Z_s\in \Dyc _q^m$, where $Q_1,\dots ,Q_r, Z_1,\dots , 
Z_s$ are prime Dyck paths.

$(1)$ For $0\leq i < j\leq m$, applying a recursive argument on $s$ and Lemma \ref{lemma2} it is easy to see that, 
\begin{align*}
P*_{\lam} (Q*_{\underline {\tau}}Z) &= ((P\t _{\lambda _1+\dots +\lambda _r}Q_1)\t _{\lambda _2+\dots +\lambda _r} \dots )\t _{\lambda _r} (Q_r*_{\underline {\tau}} Z) =\\
&(((P\t _{\lambda _1+\dots +\lambda _r}Q_1)\t _{\lambda _2+\dots +\lambda _r} \dots )\t _{\lambda _r} Q_r)*_{\underline {\delta }} Z,\end{align*}
for any pair $(\lam ,{\underline {\tau}})\in \Lambda _r^i(P)\t \Lambda _s^j(Q)$, where ${\underline {\delta }}= (\tau _0,\dots ,\tau _{s-1}, \tau_s + \lambda _r)$.

Applying the same notation than in Proposition \ref{proprelac}, we get that 

\noi $P*_{\lam} (Q*_{\underline {\tau}} Z) = (P*_{\lam} Q)*_{\underline {\delta }} Z$ if, and only if, $\psi _{ij}(P,Q) (\lam ,{\underline {\tau}}) = (\lam , {\underline {\delta}})$. 

The result follows applying point $(1)$ of Proposition \ref{proprelac}.
\medskip

$(2)$ We write ${\displaystyle \sum _{j = 0}^i P*_i (Q*_jZ) = P*_i (Q*_0Z) + \sum _{j=1}^i P*_i (Q*_j Z)}$ and 
we work the terms on the right hand side separately.
\medskip

$(a)$ Suppose that ${\underline {\tau}}\in \Lambda _s^0(Q)$, by Lemma \ref{lemma2} we get
\begin{equation*} Q*_{\underline {\tau}}Z = Q_1\t_0\dots \t _0Q_{r-1}\t_0 (Q_r*_{\underline {\tau}\rq} (Z_1\t _0\dots \t _0 Z_{j_{\underline {\tau}}}))\t_0 Z_{j_{\underline {\tau}}+1}\t _0\dots \t_0 Z_s,\end{equation*}
where $\underline {\tau}\rq = (\tau _0,\dots ,\tau _{j_{\underline {\tau}}})$ and $Q_r*_{\underline {\tau}\rq} (Z_1\t _0\dots \t _0 Z_{j_{\underline {\tau}}})$ is prime.

Applying $P*_{\lam}$, we obtain that
\begin{align*}
P&*_{\lam} (Q*_{\underline {\tau}}Z) =\\
&(P*_{\lam^1} (Q_1\t _0\dots \t _0 Q_{r-1}\t_0 (Q_r*_{\underline {\tau}\rq} (Z_1\t _0\dots \t _0 Z_{j_{\underline {\tau}}}))))*_{\lam^2} (Z_{j_{\underline{\tau}}+1}\t_0\dots \t _0 Z_s)=\\
& ((P*_{\lam ^1}Q)*_{\underline {\tau}^2}(Z_1\t_0 \dots \t _0 Z_{j_{\underline {\tau}}}))*_{\lam^2} (Z_{j_{\underline{\tau}}+1}\t_0\dots \t _0 Z_s)=(P*_{\lam ^1}Q)*_{\underline {\delta}} Z,\end{align*}
for the weak compositions ${\lam^1} = (\lambda _0,\dots ,\lambda _{r-1},\lambda _r+\dots +\lambda _{r+s{-}j_{\underline {\tau}}})$, 

\noi $\lam^2 = (\lambda _r,\dots ,\lambda  _{r+s{-}j_{\underline {\tau}}})$, ${\underline {\tau}^2}= (\tau_0,\dots ,\tau_{j_{\underline {\tau}}-1}, \tau_{j_{\underline {\tau}}}+ \lambda _r+\dots +\lambda  _{r+s{-}j_{\underline {\tau}}})$ and
${\underline {\delta}} = (\tau_0,\dots ,\tau_{j_{\underline {\tau}}-1}, \tau_{j_{\underline {\tau}}}+ \lambda _r,\lambda _{r+1},\dots ,\lambda  _{r+s{-}j_{\underline {\tau}}})$.
\medskip

The formula above implies that for any pair $(\lam , {\underline {\tau}})\in \Lambda _{r+s{-}j_{\underline {\tau}}}^i(P)\t \Lambda _s^0(Q)$, the elements $P*_{\lam}(Q*_{\underline {\tau}}Z)$ and $(P*_{\underline {\gamma}} Q)*_{\underline {\delta}}Z$ are equal whenever 
\begin{equation*}\psi _i^1(P,Q)(\lam , {\underline {\tau}}) = ({\underline {\gamma}},{\underline {\delta}}).\end{equation*}
So, we have proved that 
\begin{equation*}P*_i(Q*_0 Z) =\sum _{({\underline {\gamma}},{\underline {\delta}})} (P*_{\underline {\gamma}} Q) *_{\underline {\delta}} Z,\end{equation*}
where the sum is taken over all 
${\underline {\gamma}}\in \displaystyle\bigcup_{j=i}^m \Lambda^j_r(P)$ and 
${\underline {\delta}}\in\Lambda _s^i(P*_{\underline {\gamma}}Q)$, satisfying that $\delta _s \leq \gamma _r$.

$(b)$ Suppose that $(\lam ,{\underline {\tau}})$ belongs to $\Lambda _r^i(P)\t ({\displaystyle \bigcup _{j=1}^i\Lambda _s^j(Q)})$. We have that
\begin{equation*}Q*_{\underline {\tau}}Z = Q_1\t _0\dots \t _0 Q_{r-1}\t _0 (Q_r*_{\underline {\tau}} Z),\end{equation*}
with $Q_1,\dots ,Q_{r-1}, Q_r*_{\underline {\tau}} Z$ prime. So, computing
\begin{align*} P*_{\lam} (Q*_{\underline {\tau}} Z)=&(P*_{\lam^1}(Q_1\t_0\dots \t _0 Q_{r-1}))\t _{\lambda _r} (Q_r*_{\underline {\tau}} Z)=\\
&(P*_{\lam} Q)*_{\underline {\delta}} Z,\end{align*}
where $\lam ^1 = (\lambda _0,\dots ,\lambda _{r-2},\lambda _{r-1}+\lambda _r)$ and ${\underline {\delta}} = (\tau _0,\dots, \tau_{s-1},\tau_s+\lambda _r).$
\medskip

 Using the notation of Proposition \ref{proprelac}, we have proved the equality
\begin{equation*}P*_{\lam} (Q*_{\underline {\tau}} Z) = (P*_{\underline {\gamma}} Q)*_{\underline {\delta}} Z,\end{equation*} whenever 
 $\psi _i^2(P,Q)(\lam ,{\underline {\tau}}) = ({\underline {\gamma}}, {\underline {\delta}})$. So, we get
\begin{equation*}\sum _{j = 1}^i P*_i(Q*_j Z) = \sum _{({\underline {\gamma}},{\underline {\delta}})} (P*_{\underline {\gamma}} Q) *_{\underline {\delta}} Z,\end{equation*}
where the sum is taken over all $({\underline {\gamma}},{\underline {\delta}})\in \Lambda _r^i(P)\t \Lambda _s^i(P*_{\underline {\gamma}}Q)$ such that $\delta _s > \gamma _r$.
\medskip

Finally, adding up $(a)$ and $(b)$, we get
\begin{equation*}\sum _{j = 0}^i P*_i (Q*_j Z) = \sum _{j=i}^m (P*_jQ)*_i Z,\end{equation*}
which ends the proof.
\end{proo}
\medskip

Theorem \ref{theorem1} asserts that the graded vector space $\K [\Dyc^m]$ spanned by the set of all $m$-Dyck paths, equipped with the operations $*_i$, is a $\Dy^m$ algebra, for all $m\geq 1$.

We now turn to prove that $(\K[\Dyc^m]; *_0,\dots ,*_m)$ is in fact the free $\Dy^m$ algebra on one generator, that means $\K[\Dyc^m]$ is isomorphic to ${\mathfrak D}^m$. 
\medskip

\begin{proposition}\label{Dyckpathsobtainedfromsmallerones} Any element of $P\in \Dyc_n^m$ is of the form $R_1*_i R_2$, where $0\leq i\leq m$ and the sizes of $R _1$ and $R_2$ are strictly smaller than 
$n$, for $n\geq 1$.
\end{proposition}
\medskip

\begin{proo} Suppose that $P=P_1\t_0 P_2\t_0\dots \t_0 P_r$, with $P_i$ prime, for $1\leq i\leq r$.

If $r>1$, then $P= P\rq *_0 P_r$, with $P\rq :=P_1\t_0\dots \t_0 P_{r-1}$, and the result is true.

If $P$ is prime, then $0\leq L_1(P) <m$. Let $X(P)=\{1< s_1<\dots < s_k=n\}$ be the set of integers satisfying that there exists at least one down steps in $\dw^{s_j}(P)$ of color $1$, and let $h_j>0$ be the number of down steps of color $1$ in 
$\dw^{s_j}(P)$, for $1\leq j\leq k$. We have that $h_1+\dots +h_k+ L_1(P) = m$.

It is immediate to see that there exist $m$-Dyck paths $Q_1,\dots ,Q_k$, satisfying that
\begin{equation*} P=(((\rho _m\t_{l_1} Q_1)\t_{l_2} Q_2)\dots )\t_{l_k}Q_k,\end{equation*}
where $l_j=h_j+\dots +h_k$, for $1\leq j\leq k$.

Let $R_1= (((\rho _m\t_{l_1} Q_1)\t_{l_2} Q_2)\dots )\t_{l_{k-1}} Q_{k-1}$. The size of $R_1$ is $n_1$, for some $0<n_1 <n$, and the number of steps of color $1$ in $\dw^{n_1}(R_1)$ is $h_{k-1}+h_k$.
As $h_{k-1}\geq 1$, we have that  $P=R_1*_{h_k}Q_k$, which ends the proof.\end{proo}
\medskip

The following theorem states that the graded vector space $\D_m$ also describes the algebraic operad $\Dy ^m$.

\begin{theorem}\label{theorem: Dm is free on one generator}
The free $\Dy^m$ algebra on one generator is isomorphic to $(\K[\Dyc^m],*_0,\dots,*_m)$.
\end{theorem}
\medskip

\begin{proo}
As $\K[\Dyc^m]$ is a $\Dy ^m$ algebra, there exists a unique homomorphism $\phi : {\mathfrak D}^m\longrightarrow \K[\Dyc^m]$ such that $\phi (\vert)$ is $\rho _m$, the unique $m$-Dyck path of size $1$.
Proposition \ref{Dyckpathsobtainedfromsmallerones} implies that $\phi$ is surjective. 

By Proposition \ref{Dyckpathsobtainedfromsmallerones}, the subspace of homogeneous elements of degree $n$ of $\K[\Dyc^m]$ is generated by the subset $\Dy_n^m$ of $m$-Dyck paths of size $n$. Let ${\mathfrak D}_n^m$ be the subspace of elements of degree $n$ of ${\mathfrak D}^m$.

As $\phi $ is surjective, to prove that $\phi$ is an isomorphism it suffices to show that the dimension of the vector space 
${\mathfrak D}_n^m$ is the number of elements of the set $\Dyc _n^m$, that is 
\begin{equation*}{\mbox {dim}_{\K} ({\mathfrak D}_n^m)}= \vert \Dyc_n^m\vert = d_{m,n},\end{equation*}
which was proved in Corollary \ref{dimensionofthefree}.
\end{proo}

\subsection{Connection with the $m$-Tamari lattice}
\medskip

F. Bergeron extended the Tamari order to the sets $\Dyc_n^m$ of Dyck paths (see \cite{BerPre}) . Let us describe briefly the $m$-Tamari lattice $\Dyc_n^m$.
\medskip

Let $P$ be an $m$-Dyck path. 
For any down step $d_0\in \dw (P)$ which is followed by an up step $u\in \up (P)$, consider the excursion $P_u$ of $u$ in $P$ and its matching down step $w_u$ as described in Definition \ref{excursion}. Let $P_{(d_0)}$ be the Dyck path obtained by removing $d_0$
and gluing the initial vertex of $u$ to the end of the step preceding $d_0$, and attaching $d_0$ at the final point of $w_u$. For example
\begin{figure}[h]
\includegraphics[scale=0.6]{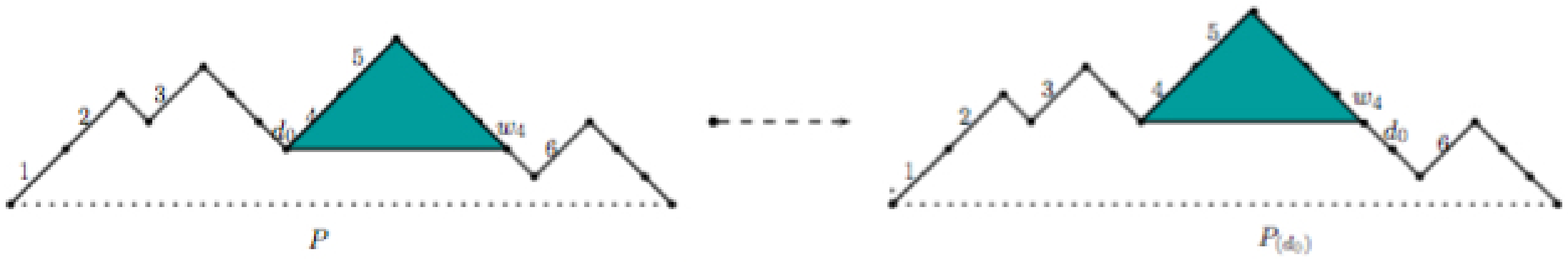}
\end{figure}

It is immediate to see that $\alpha _{P_{d_0}}(d) = \alpha _P(d)$, for any $d\in \dw (P)$.

\begin{definition} \label{defmTamari}  The $m$-Tamari order on $\Dyc_n^m$ is the transitive relation spanned by the covering relation
\begin{equation}P \lessdot P_{(d)},\end{equation}
for any $d\in \dw (P)$ such that the final vertex of $d$ is the initial point of an up step $u\in \up (P)$. We use the symbol $\lessdot $ for a covering relation. \end{definition}

The Hasse diagrams for $m = 2$ and $n = 1,2$ are

\begin{figure}[h]
\includegraphics[scale=0.6]{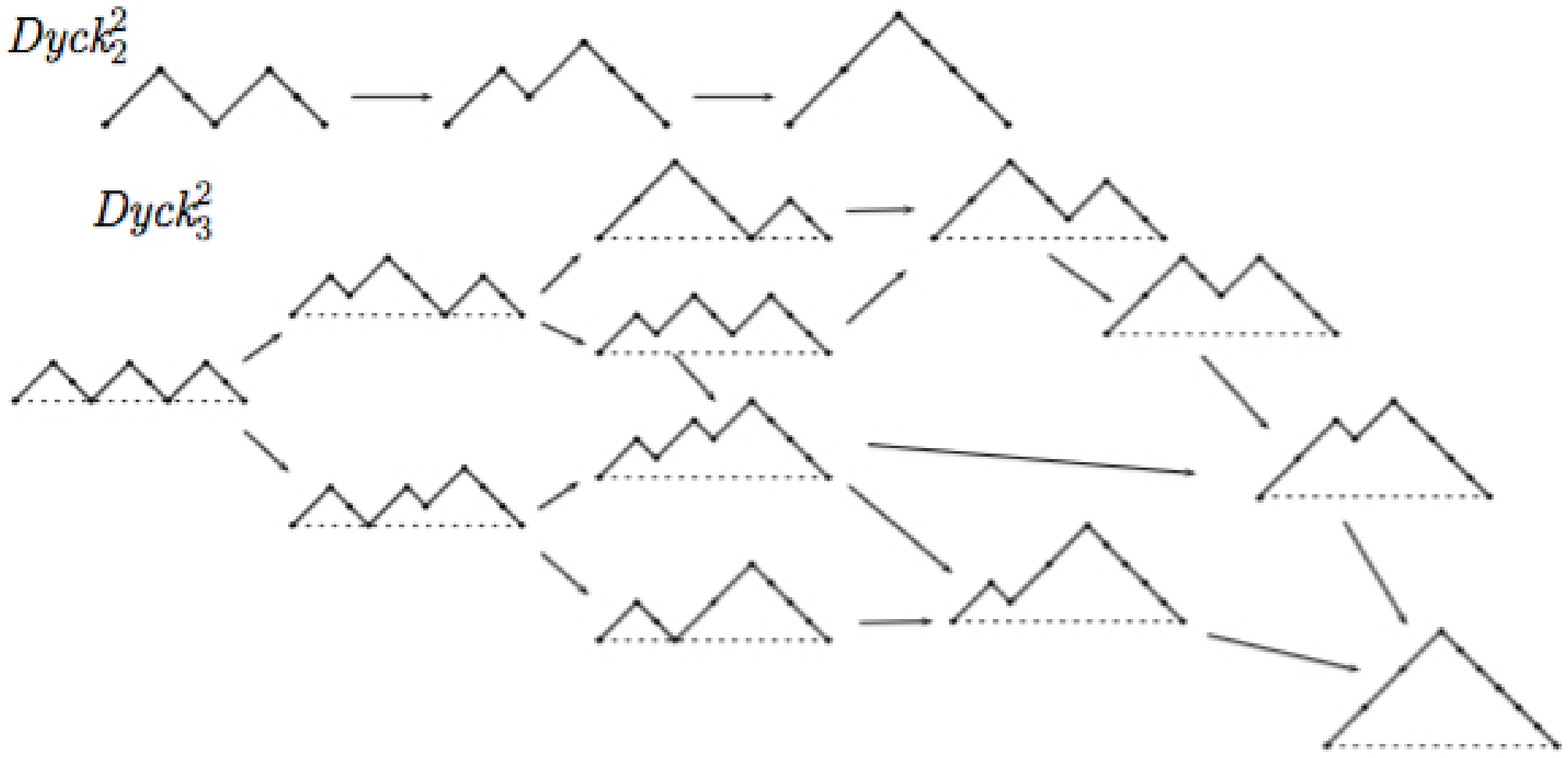}
\end{figure}
\medskip

The goal of the present section is to show that the binary operations $*_i: \K[\Dyc_n^m]\ot \K[\Dyc_r^m]\longrightarrow \K[\Dyc_{n+r}^m]$ 
are described in terms of the $m$-Tamari order.

\begin{remark} \label{remmatching} \begin{enumerate} \item Let $Q$ be a prime Dyck path, for any pair of $m$-Dyck path $P$, we get
\begin{equation} P\t _0 Q < P\t _1 Q <\dots < P\t _{L(P)}Q,\end{equation}
in the $m$-Tamari lattice.
\item If $P < P\rq $ in $\Dyc_{n_1}^m$ are such that $L(P) = L(P\rq)$, and $Q < Q\rq $ in $\Dyc_{n_2}^m$, then \begin{enumerate}
\item $P\t _k Q < P\t _k Q\rq $, for any $0\leq k\leq L(P)$, 
\item $P\t _k Q < P\rq \t_k Q$, for any $0\leq k\leq L(P)$.\end{enumerate}
\end{enumerate}
\end{remark}

For the rest of the section, the $m$-Dyck path $Q$ is supposed to be a product $Q=Q_0\t _0\dots \t_0 Q_r$, where all the $Q_j$\rq s are prime Dyck paths.

\begin{lemma} \label{lemmalambdaordre} Let $P\in \Dyc_{n_1}^m$ and $Q\in \Dyc_{n_2}^m$ be two Dyck paths. Two weak compositions $\lam$ and ${\underline {\gamma}}$ in $\Lambda_r(P)$ satisfy that 
\begin{equation*}\lambda _j + \dots +\lambda _r \leq \gamma _j + \dots + \gamma _r,\end{equation*}
for $1\leq j\leq r$, if, and only if, $P*_{\lam}Q \leq P*_{\underline {\gamma}} Q$.\end{lemma}
\medskip

\begin{proo} If $Q$ is prime, the result follows from point $(1)$ of Remark \ref{remmatching}. 
Suppose that $Q = Q_1\t _0\dots \t_0 Q_r$, for $r >1$. 
A recursive argument shows that, for any pair of elements $\lam\rq$ and ${\underline {\gamma }\rq}$ in $\Lambda _{r-1}(P)$,  we have that 
\begin{equation*}P*_{\lam\rq} (Q_1\t_0 \dots \t_0 Q_{r-1}) \leq P*_{\underline {\gamma}\rq} (Q_1\t_0 \dots \t_0 Q_{r-1}),\end{equation*}
 whenever $\lambda\rq _j + \dots +\lambda\rq  _{r-1} \leq \gamma\rq  _j + \dots + \lambda\rq _{r-1},$ for $1\leq j\leq r-1$.
\medskip 
 
So, we get\begin{enumerate}[(a)]
\item $P*_{\lam}Q = (P*_{\lam\rq} (Q_1\t _0\dots \t_0 Q_{r-1}))\t_{\lambda _r} Q_r$, 
\item $P*_{\underline {\gamma}} Q = (P*_{\underline {\gamma}\rq} (Q_1\t _0\dots \t_0 Q_{r-1}))\t_{\gamma _r} Q_r$, 
\end{enumerate}
where $\lam \rq = (\lambda _0, \dots ,\lambda _{r-1}, \lambda _{r-1} + \lambda _r)$ and 
$\underline {\gamma}\rq = (\gamma _0, \dots ,\gamma _{r-1}, \gamma _{r-1} + \gamma_r)$.
The recursive hypothesis implies that 
\begin{equation*}P*_{\lam\rq} (Q_1\t _0\dots \t_0 Q_{r-1})\leq P*_{\underline {\gamma}\rq} (Q_1\t _0\dots \t_0 Q_{r-1}),\end{equation*}
 and, using that $\lambda _r \leq \gamma _r$, we obtain $P*_{\lam}Q \leq P*_{\underline {\gamma}} Q$.

Conversely, suppose that $P*_{\lam}Q \leq P*_{\underline {\gamma}} Q$. Point $(3)$ of Remark \ref{remmatching} implies that 
\begin{equation*}\lambda _j + \dots +\lambda _r \leq \gamma _j + \dots + \lambda _r,\end{equation*}
for $1\leq j\leq r$, which ends the proof. \end{proo}

\begin{notation} \label{notTamari} For any $m$-Dyck path $P$ of size $n$ and any $0\leq i\leq m$, let \begin{enumerate}
\item $c_i(P)$ be the minimal number of elements such that the word 
\begin{equation*}\alpha _P(d_{L(P){-}c_i(P)+1})\dots \alpha _P(d_{L(P)})\end{equation*} contains $i$ times an integer in $\{1,\dots ,n\}$ and no integer more than $i$ times,
\item $C_i(P)$ be the maximal integer such that the word \begin{equation*}\alpha _P(d_{L(P){-}C_i(P)+1})\dots \alpha _P(d_{L(P)})\end{equation*} contains at least one integer repeated $i$ times and no integer repeated $i+1$ times.\end{enumerate}

Let $P\in \Dyc_{n_1}^m$ and $Q\in \Dyc_{n_2}^m$ be two Dyck paths. For any integer $0\leq i\leq m$, let $P /_i Q$ and $P \backslash _i Q$ be the Dyck paths defined as follows\begin{enumerate}[(a)]
\item $P/_i Q := P\t _{c_i(P)} Q$,
\item $P\backslash _i Q := (P\t _{L(P)} (Q_1\t_0\dots \t _0 Q_{r-1}))\t _{C_i(P)} Q_r$.\end{enumerate}
\end{notation}

\begin{proposition} \label{propTamarim} For any pair of Dyck paths $P\in \Dyc_{n_1}^m$ and $Q\in \Dyc_{n_2}^m$ and any integer $0\leq i\leq m$, the product $*_i$ is given in terms of the $m$-Tamari order by the following formula
\begin{equation}P*_iQ = \sum _{P/_iQ\leq Z\leq P\backslash _i Q} Z.\end{equation}\end{proposition}
\medskip

\begin{proo} Suppose that $Q=Q_1\t_0\dots \t_0Q_r$, with all the $Q_i$\rq s prime and that $\lam \in \Lambda _r^i(P)$. 

The weak composition $\lam = (\lambda _0,\dots ,\lambda _r)$ satisfies that $c_i(P)\leq \lambda _r\leq C_i(P)$ and ${\displaystyle \sum _{j=0}^r\lambda _i = L(P)}$. 
As we have that\begin{enumerate}[(a)]
\item $P/_iQ = P*_{(L(P){-}c_i(P), 0, \dots ,0, c_i(P))} Q$, and 
\item $P\backslash _i Q = P*_{(0,\dots , 0, L(P){-}C_i(P), C_i(P))}Q$,\end{enumerate}
applying Lemma \ref{lemmalambdaordre}, it is easily seen that $P/_i Q\leq P*_{\lam}Q\leq P\backslash _i Q$.

Recall that, whenever $R < S$ in the Tamari lattice, the set $\dw (R)$ of down steps of $R$ is identified with the set $\dw  (S)$. For any $d\in \dw (P)$ the levels of $d$ in $R$ and in $S$ are different but $\alpha _R(d) = \alpha _S(d)$. 

Note that the unique down steps which have different levels in the Dyck paths $P/_iQ$ and $P\backslash _i Q$ are colored by the set of integers $\{ 1,\dots ,n_1\}$. So, for any $P/_iQ\leq Z\leq P\backslash _i Q$ and any $1\leq l\leq r$, we get that 
\begin{equation}L_j(Z) =L_j(Q_l),\ {\rm for}\ n_1+n_{21}+\dots +n_{2(l-1)} < j < n_1+n_{21}+\dots + n_{2l}.\end{equation}

Define \begin{equation*}\lambda _j = \begin{cases} L_{n_1+n_{21}+\dots +n_{2j}} (Z) -L(Q_j),&\ {\rm for}\ 1\leq j\leq r,\\
L_{n_1}(Z) - L(P),&\ {\rm for}\ j=0.\end{cases} \end{equation*}

The arguments above show that \begin{enumerate}
\item $c_i\leq \lambda _r\leq C_i$,
\item $c_i \leq  \lambda _j+ \dots +\lambda _r \leq L(P)$, for $1\leq j\leq r-1$,
\item $0\leq L_{n_1}(Z)\leq L(P) - c_i.$\end{enumerate}

From $(4.4.4)$, we get that $Z = P*_{\lam} Q$.

Lemma \ref{lemmalambdaordre} and $P/_iQ\leq P*_{\lam}Q\leq P\backslash _i Q$ imply that $\lam \in \Lambda _r^i(P)$.\end{proo}

Let us define the product $*$ on $\K[\Dyc^m]$ as the sum $* := *_0+\dots +*_m$. It is not difficult to see, using Proposition \ref{propTamarim}, that 
\begin{equation*}P*Q = \sum _{P/_0 Q \leq Z\leq P\backslash _mQ} Z.\end{equation*}

\begin{example} Consider the Dyck paths $P = (1,3)$ and $Q=(2,2)$ in $\Dyc_2^2$, the following diagram describes the Tamari interval $I_{P*Q}$ of all $Z\in \Dyc_4^2$ such that 
$P * Q = {\displaystyle \sum _{Z\in I_{P*Q}} Z}$.

The Dyck paths in red are the terms of $P*_0Q$, the ones in green are the terms of $P*_1Q$, 
and the ones in blue are the terms appearing in $P*_2Q$.
\medskip

\begin{figure}[h]
\includegraphics[scale=0.6]{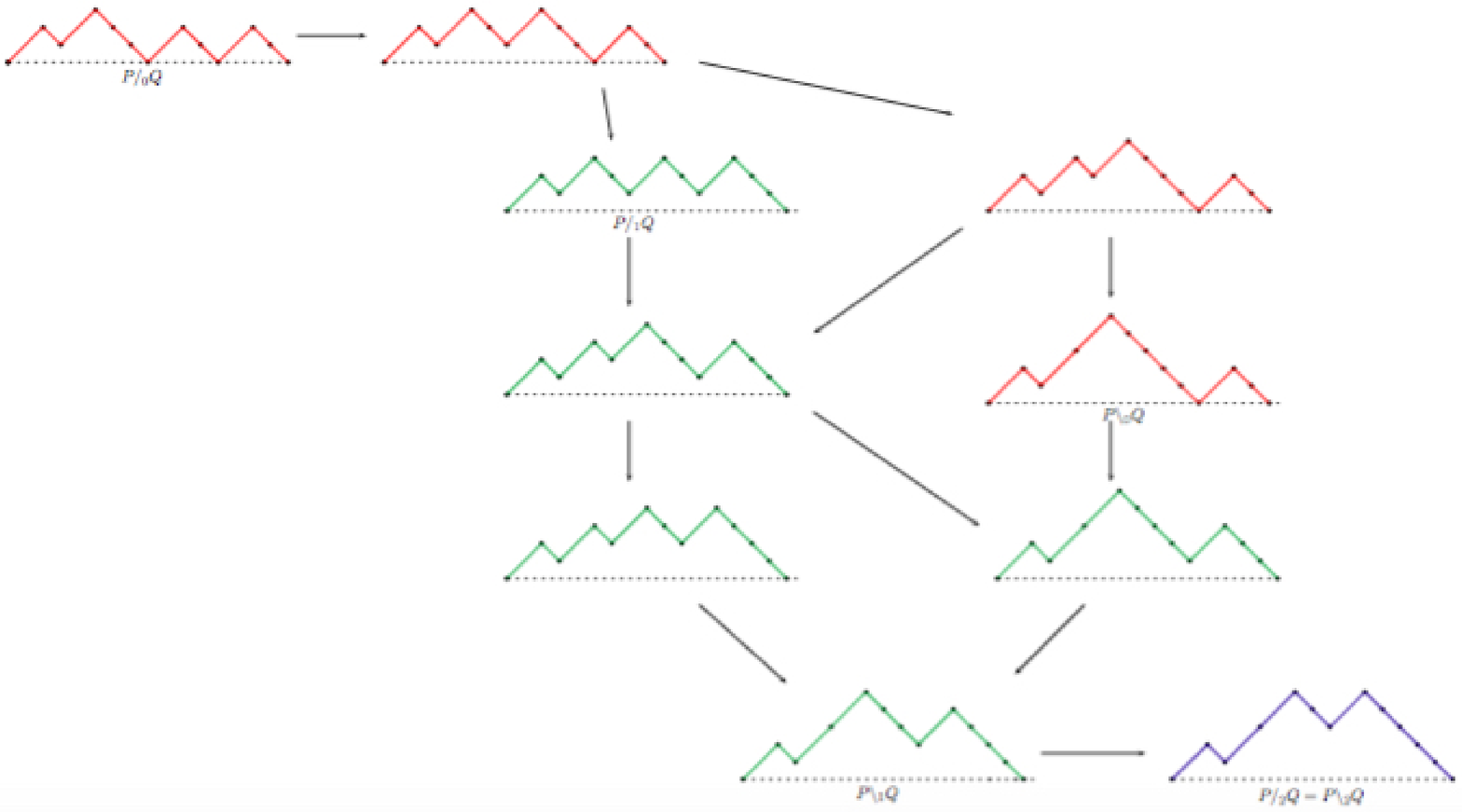}
\end{figure}
\end{example}

\end{document}